%%%%%%%%%%%%%%%%%%%%%%%%%%%%%%%%%%%%%%%%%%%%%%%%    
%
%        THIS IS A  PLAIN TeX FILE
%
%%%%%%%%%%%%%%%%%%%%%%%%%%%%%%%%%%%%%%%%%%%%%%%%

\magnification=1200

\font\titfont=cmr10 at 12 pt

\font\headfont=cmr10 at 12 pt

%\font\AAA=Times.dfont  at 12pt
 %\font\BBB=Times.dfont at 8pt

%\font\AAA=cmr10 at 12pt
%\font\BBB=cmr10 at 8pt

\overfullrule=0in

\def\boxit#1{\hbox{\vrule
 \vtop{%
  \vbox{\hrule\kern 2pt %
     \hbox{\kern 2pt #1\kern 2pt}}%
   \kern 2pt \hrule }%
  \vrule}}
  \def\mathqed{  \vrule width5pt height5pt depth0pt}

  \def\harr#1#2{\ \smash{\mathop{\hbox to .3in{\rightarrowfill}}\limits^{\scriptstyle#1}_{\scriptstyle#2}}\ }

\def\bra#1#2{\langle #1, #2\rangle}
\def\bbf{{\bf F}}

\def\ss{\subset}

\def\half{\hbox{${1\over 2}$}}
\def\smfrac#1#2{\hbox{${#1\over #2}$}}

\def\log{{\rm log}}

\def\tr{{\rm tr}}
\def\max{{\rm max}}
\def\min{{\rm min}}

\def\det{{\rm det}}

\def\Sym{{\rm Sym}^2}

\def\rn{\bbr^n}

\def\Int{{\rm Int}}

\def\Symn{{\Sym(\rn)}}

\def\Theorem#1{\medskip\noindent {\bf THEOREM \bf #1.}}
\def\Prop#1{\medskip\noindent {\bf Proposition #1.}}
\def\Cor#1{\medskip\noindent {\bf Corollary #1.}}
\def\Lemma#1{\medskip\noindent {\bf Lemma #1.}}
\def\Remark#1{\medskip\noindent {\bf Remark #1.}}

\def\Def#1{\medskip\noindent {\bf Definition #1.}}

\def\Ex#1{\medskip\noindent {\bf Example \bf    #1.}}

\def\pf{\medskip\noindent {\bf Proof.}\ }
\def\qed{\hfill  $\vrule width5pt height5pt depth0pt$}

   \def\cp{{\cal P}}

\def\cp{{\cal P}}

\def\vf{\varphi}

\def\wt{\widetilde}

\def\and{\qquad {\rm and} \qquad}

\def\bbr{{\bf R}}\def\bbh{{\bf H}}
\def\bbc{{\bf C}}

\def\a{\alpha}
\def\b{\beta}
\def\d{\delta}
\def\e{\epsilon}

\def\l{\lambda}

\def\s{\sigma}

\def\L{\Lambda}

\def\O{\Omega}

\def\lloc{L^1_{\rm loc}}

\def\bo{\partial \Omega}

\def\Symn{\Sym(\rn)}
 
\def\USC{{\rm USC}}

\def\cpt{\wt{\cp}}
\def\ft{\wt F}
\def\ob{\overline{\O}}

\def\AA{1}
\def\BB{2}
\def\CC{3}

\def\EE{4}
\def\FF{5}
\def\GG{5}
\def\HH{6}

\def\II{7}
\def\JJ{8}
 \def\KK{9}
 \def\LL{10}
\def\MM{11}

 \def\oL{{\overline \L}}
 \def\uL{{\underline \L}}
\def\of{\overline f} 
\def\uf{\underline f}

  \def\AS{1}
  \def\BD{2}
  \def\BaB{3}
  \def\BRUCK{4}
  \def\CaC{5}
  \def\CLN{6}
  \def\C{7}
  \def\CIL{8}
  \def\EB{9}
  \def\F{10}
  \def\DDD{11} 
  \def\DDR{12}
 \def\SURVEY{13}
 \def\RS{14}
 \def\ASPECTS{15}
\def\AETHM{16}
\def\HIR{17}
 \def\Kawohl{18}
 \def\Krylov{19}
 \def\LE{20}

%\centerline{\titfont  THE MAXIMUM PRINCIPLE  }
%\medskip

\centerline{\titfont   CHARACTERIZING THE STRONG MAXIMUM PRINCIPLE}
\medskip

\centerline{\titfont   FOR CONSTANT COEFFICIENT SUBEQUATIONS}
\vskip .2in

\centerline{\titfont F. Reese Harvey and H. Blaine Lawson, Jr.$^*$}
\vglue .9cm
\smallbreak\footnote{}{ $ {} \sp{ *}{\rm Partially}$  supported by
the N.S.F. } 
%\vskip .3in

\centerline{\bf ABSTRACT} \medskip
  \font\abstractfont=cmr10 at 10 pt
%For subsolutions of a degenerate elliptic equation $\bbf(D^2)=0$,
%several (equivalent) characterizations of the cases where the
%maximum principle holds were given in [DD].  When the maximum 
%principle holds, the strong maximum principle may or may not hold.  Here 
%we introduce an easily computed function $f$   on $(0, \infty)$ which is
%defined by $\bbf$, and we show that the strong maximum principle holds 
%depending on whether $\int_{0^+} {dy \over f(y)}$ is infinite or finite.

{{\parindent= .83in
\narrower\abstractfont \noindent
In this paper we characterize the degenerate elliptic equations $\bbf(D^2u)=0$
whose subsolutions ($\bbf(D^2u)\geq 0$) satisfy
the strong maximum principle.  We introduce an easily computed function
$f$ on $(0, \infty)$ which is determined by $\bbf$, and we show that the 
strong maximum principle holds 
depending on whether $\int_{0^+} {dy \over f(y)}$ is infinite or finite.
This is in the spirit of previous work characterizing the ordinary maximum
principle in terms of the geometry of the set of symmetric matrices
$F = \{\bbf \geq 0\}$. Along the way, radial subsolutions are characterized,
and,  as an application, a sufficient condition for strong  comparison is  established.
A number of examples, important for the theory of such equations, are examined.

}}

\vskip.4in

\centerline{\bf TABLE OF CONTENTS} \bigskip

{{\parindent= .1in\narrower\abstractfont \noindent

\qquad \AA. Introduction.\smallskip

\qquad \BB.    Characterizing the Maximum Principle.    \smallskip

\qquad \CC.     Characterizing the Strong Maximum Principle -- Three Cases.    \smallskip

%\qquad \DD.  Characterizing the (SMP) for Cone Subequations.
%     \smallskip

\qquad \EE.    The Radial Subequation Associated to $F$.    \smallskip

%\qquad \FF.   Radial Subharmonics.    \smallskip

\qquad \GG.  Increasing Radial Subharmonics for Borderline Subequations.     \smallskip

\qquad \HH.   Proof of the (SMP) in the Borderline Case.     \smallskip

\qquad \II.   Radial (Harmonic)  Counterexamples to the (SMP).    \smallskip

\qquad \JJ.  Subequations with the Same Increasing Radial Subharmonics.    \smallskip

\qquad \KK.  Strong Comparison and Monotonicity.     \smallskip

\qquad\LL. Examples of Exotic Monotonicity Subequations which are not Cones.\smallskip
   
\qquad \MM.\  Another Application  --  Product Subequations.
  
}}

\vskip .15in

\qquad \qquad
{\bf Appendix A.} A Theorem on Radial Subharmonics.

\smallskip
\qquad \qquad
{\bf Appendix B.} Uniform Ellipticity and the Borderline Condition.

\vfill\eject

\noindent{\headfont \AA.\  Introduction}
\medskip

This paper is concerned with differential equations of the form $\bbf(D^2u) =0$
where $\bbf$ is degenerate elliptic, and  attention is focused on the set
$F(X)$ of viscosity subsolutions ($\bbf(D^2u) \geq 0$) on an open set $X\ss\rn$
as defined in the seminal papers [\C] or [\CIL].  
Our interest here is in the theory of such equations, rather than concern with specific cases.
%The scope of the paper is quite limited in the sense that $\bbf$ is constant
%coefficient, pure second-order.  However, it is broad in 
%that the usual structural assumptions
%such as uniform ellipticity  or homogeneity are not required.

The main point of the paper
is to address the following.
\medskip
\noindent
{\bf Question:} 
When do the subsolutions satisfy the strong maximum principle?\medskip

 By the {\bf maximum principle}
 and the {\bf strong maximum principle} for $\bbf$ we mean the following.
 Given a bounded domain $\O\ss\rn$, let $F(\ob)$ denote the space of 
 upper semi-continuous functions  on $\ob$ which are subsolutions on $\O$.
 Consider the implications:
 $$
 u\in F(\ob)\qquad\Rightarrow \qquad \sup_{\ob} u \ \leq\ \sup_{\bo} u
 \eqno{(MP)}
 $$
 $$
 u\in F(\ob) \ \ {\rm has\ an\ interior\ maximum\ point} \qquad\Rightarrow \qquad u\ \ {\rm is\ constant}
 \eqno{(SMP)}
 $$
We say that the (MP)/(SMP) holds for $F$ if it is true for all such $\O$ and $u$.
Of course, (SMP) $\Rightarrow$ (MP).

\medskip

There is, of course, a huge literature concerned with the (SMP) for viscosity subsolutions
of nonlinear equations (for just a few examples see  [\BD], [\BaB], [\CaC], [\Kawohl], [\LE] and the references therein.
We note in particular the landmark paper [\BaB] and refer the reader to Remark \JJ.10 for
a discussion of its relationship to the results here.)
However, these authors typically make  structural assumptions
such as uniform ellipticity  or homogeneity.

Here we confine our attention to the special  case of 
constant coefficient, pure second-order (degenerate elliptic) equations in $\rn$,
sometimes having an invariance property.
No other structural assumptions are made.
The point of this paper is that in this special but important setting one can give
a complete and somewhat unexpected   characterization of exactly when the 
(SMP) holds.

This work is a natural outgrowth of the results in [\DDD] where the ordinary (MP)
is characterized, in several equivalent ways, in terms of the geometry of the set
$F\equiv \{A : \bbf(A)\geq0\}$. They are given   in Theorem \BB.1 below.
We have reviewed these  (MP) results   in Section \BB\  for two reasons.  
First, they are scattered  in various remarks throughout [\DDD]. 
Secondly, these geometrically  elegant and complete characterizations 
available for the (MP) provide the motivation for our discussion of the (SMP) in this paper.

%A   characterization of when the (MP) holds for $\bbf$ was given in [\DDD, Remark 4.7].
%This theorem is amplified in Section \BB.

Our discussion of the (SMP) divides into three cases.  The first case is relatively simple and classical,
and the (SMP) always holds:
$$
\bbf(0)\ <\ 0\qquad\Rightarrow\qquad {\rm The\ (SMP)\ holds.}
\eqno{(\AA.1)}
$$

The second case is also simple and classical, and the (SMP) always fails.  To describe this case
we set some notation.   Given a non-zero 
vector $e\in\rn$, let $P_e$ and $P_{e^\perp}$ denote orthogonal projection onto the 
line spanned by $e$ and the hyperplane perpendicular to $e$ respectively, so that
$P_e+ P_{e^\perp} = I$. Our second case is the following.
$$
\bbf(-\mu P_e) \geq 0 \ \ {\rm for\ some \ } \mu >0, e\neq0
 \quad \Rightarrow\quad {\rm (SMP)\ fails}.
 \eqno{(\AA.2)}
$$

Consequently,  we concentrate on the remaining case, which will be referred to as the {\bf borderline case}. 

In this paper a key role is played by the increasing radial subsolutions.
They are determined by a ``characteristic function'' $f$ of $\bbf$, 
which is defined as follows.   For simplicity we first assume  the following weak form of
{\bf invariance}: For all $\l,\mu\in\bbr$,
$$
\eqalign{
&\bbf(\l P_{e^\perp} -\mu P_e)\ \geq\ 0 \quad {\rm for\ one\ }\ e\ne0  \cr
\Rightarrow \qquad &\bbf(\l P_{e^\perp} -\mu P_e)\ \geq\ 0 \quad {\rm for\ all\ \  }\ e\ne0 
}
 \eqno{(\AA.3)}
$$
This holds, for example,  if $\bbf$ is invariant under a group such as O$_n$ or SU$_{n/2}$
acting transitively  on the $n-1$ sphere in $\rn$ (a condition called {\sl ST-invariance} in [\ASPECTS]).
Given an invariant  $\bbf$,   the {\bf characteristic function}  $f$ associated to $\bbf$  for $0\leq \l<\infty$ 
is defined by
$$
f(\l) \ \equiv\ \sup\left\{  \mu : \bbf\left (\l P_{e^\perp} -\mu P_e\right)\ \geq\ 0     \right\}
\eqno{(\AA.4)}
$$
The {\bf borderline}  cases are exactly the cases where $f(0)=0$ (see Lemma \CC.4).

Now we can state  our main result, simplified by assuming invariance.

\Theorem {A}  {\sl
Suppose $\bbf$ is invariant  and borderline. Then}
$$
{\sl The \rm \  (SMP)\ \sl holds\ for\ } \bbf 
\qquad\iff\qquad 
\int_{0^+} {dy \over f(y)} \ =\ \infty.
$$

The general (non-invariant) version of this result is given below in Theorem A$'$.

The characteristic  function $f$ determines the following one-dimensional variable coefficient 
operator
$$
\left(R^{\uparrow}_f \psi\right) (t) \ \equiv \ \min\left\{   \psi'(t), \ \psi''(t)+f\left({\psi'(t)\over t}\right)\right\}.
\eqno{(\AA.5)}
$$
The next result is of general interest, and probably classical in the $C^2$-case.

\Prop{B} {\sl
A radial function $u(x) = \psi(|x|)$ with $\psi$   increasing,
is an $\bbf$-subsolution if and only if $\psi$ is an $R^{\uparrow}_f $-subsolution.}

\medskip
The ``only if'' part of this result requires a technical lemma for general 
upper semi-continuous functons, which is given in  Appendix A.

These two results lead to the following.
\medskip
\noindent
{\bf Question:}  Given an  upper semi-continuous, increasing function  $f: [0,\infty) \to [0,\infty]$
with $f(0)=0$, is there a way to describe all the equations $\bbf$ which have $f$ as their characteristic function,
or equivalently (by Proposition B) have the same set of increasing radial subsolutions.
\medskip

This question is answered here.   First, such equations $\bbf$ always exist for any 
such $f$ (as above).
%increasing upper semi-continuous function $f$ with $f(0)=0$.
Here are two crucial examples.  
Let $\l_1(A)\leq\cdots\leq \l_n(A)$ denote the ordered eigenvalues of a symmetric matrix $A$
so that $\l_{\rm min} = \l_1$ and $\l_{\rm max} = \l_n$.
Define
$$
\bbf^{\rm min/max}_f \left( D^2u \right) \ \equiv\
 \min \left\{\l_{\rm max}  \left( D^2u \right), \l_{\rm min}   \left( D^2u \right)
  + f\left(\l_{\rm max}  \left( D^2u \right)\right)   \right\},\ \ 
 {\rm and}
$$
$$
\bbf^{\rm min/2}_f \left( D^2u \right) \ \equiv\
 \min \left\{ \l_2 \left( D^2u \right), \l_{\rm min}  \left( D^2u \right) + f\left(\l_2 \left( D^2u \right)\right)   \right\}.
$$
Both have associated characteristic function $f$ (see Lemma \JJ.3). 
 In fact, they are the largest and the smallest such examples.
 (A priori it is not clear that there is a largest or a smallest.)
 
 \Theorem{C} {\sl 
 If $\bbf$ is invariant and borderline with characteristic function $f$, then its subsolutions satisfy
  $$
F^{\rm min/2}_f  (X)\ \ss\ F(X)\ \ss\ F^{\rm min/max}_f (X).
 $$
Conversely, these containments imply that $\bbf$ must have  characteristic function $f$.}

\medskip

Our first main result, Theorem A above, extends to  $\bbf$'s which are not necessarily invariant
 as follows.
 We define the 
 {\sl upper    and lower    characteristic functions $\of$ and $\uf$
for $\bbf$} by:
$$
\of(\l) \ \equiv \ \sup\left\{ \mu :  \bbf\left(\l P_{e^\perp} - \mu P_{e} \right)\geq0 \ \ {\rm for\ some\ } e\neq 0\right\}
$$
$$
\uf(\l) \ \equiv \ \sup\left\{ \mu :  \bbf\left(\l P_{e^\perp} - \mu P_{e} \right)\geq0 \ \ {\rm for\ all\ } e\neq 0\right\}
$$
When $\bbf(0)=0$, we have
$
  \uf(0) \ =\ \of(0) \ =\ 0
$.

\Theorem{A$'$} {\sl Suppose that $F$ is borderline  and has upper and lower characteristic
functions $\of$ and $\uf$.
$$
{\rm (a)} \quad {\sl If}\ \ \int_{0^+} {dy\over \of(y)}\ =\   \infty, \ \ {\sl then\  the \ (SMP)\ holds\ for\ } \bbf. 
$$
$$
{\rm (b)} \quad {\sl If}\ \ \int_{0^+} {dy\over \uf(y)}\ <\  \infty, \ \ {\sl then\  the \ (SMP)\ fails\ for\ }  \bbf.
 $$
}
\medskip

Now that our main result has been stated in the traditional manner using nonlinear 
operators $\bbf$, we switch to the  geometric point of view (pioneered by  Krylov [\Krylov]) 
which replaces $\bbf$ with the 
subset $F=\{\bbf \geq 0\}$ of $\Symn$, the space of $n\times n$ symmetric matrices.
This is particularly appropriate for discussing questions such as ours concerning
the (MP) and (SMP) since they only depend on the space of 
subsolutions $F(X)$ which in turn only depends on the geometry of the subset $F$ and not on the operator
$\bbf$ used to define it.  (The situation is analogous to studying submanifolds independently
of any implicit defining function.)
Let 
$$
\cp \ \equiv\  \{A\in \Symn : A\geq0\}.
$$
  Instead of ``operators'' 
we consider {\bf subequations}  which by definition are closed subsets $F\ss \Symn$
satisfying the weakest form of ellipticity, namely:
$$
F+\cp \ \ss\ F,
\eqno{(P)}
$$
called {\bf positivity}.  Subsolutions are defined in the usual manner, except
that one requires $D^2_x\vf \in F$, rather than $\bbf(D^2_x \vf)\geq0$,
for test functions $\vf$ at $x$.  To emphasize the parallels with potential
theory in several complex variables, we will use the terminology
{\bf $F$-subharmonic} rather than {\sl $\bbf$-subsolution}.
The key topological property of $F$ is that:
$$
F\ =\ \overline{\Int F}.
\eqno{(T)}
$$
This follows easily from (P) and the assumption that $F$ is closed.

\medskip
\noindent
{\sl Some Technical Points:} With the operator $\bbf$ replaced by the set
$F\equiv \{\bbf\geq0\}$, the positivity condition for $F$ is weaker than degenerate
ellipticity for the operator $\bbf$.  Positivity  is equivalent to requiring that: 
$\bbf(A)\geq0 \ \ \Rightarrow\ \ \bbf(A+P)\geq0$ for all  $A\in \Symn, P\in \cp$.
(Weak  ellipticity is the requirement that $\bbf(A+P) \geq \bbf(A)$ for all such $A$ and $P$.)

 Our notion of a supersolution $v$ is (for some $\bbf$) more restrictive
than the classical notion $\bbf(D^2 v)\leq 0$.  
We require $-v$ to be subharmonic for the dual subequation $\ft = -(\sim\Int F)$.
This has an advantage over the standard notion of supersolution.  For
example, we were able to prove that comparison always holds for 
any subequation $F\ss \Symn$ (Theorem 6.4 in [\DDD]).
For degenerate elliptic operators $\bbf$, the statement becomes: 
comparison holds if and only if $\{\bbf\leq0\}$ is the complement of the interior of $\{\bbf\geq0\}$.
The reader is referred to the ``Pocket Dictionary'' in Appendix A of  [\SURVEY]
for a more complete translation of concepts.

By {\sl  comparison} holding for a subequation $F$ we mean
$$
u+w \ \ {\rm satisfies \ the\ (MP)\ for\ all\ }\ u\in F(\ob), \ w\in \ft(\ob).
\eqno{(C)}
$$
By  {\sl strong comparison} holding for a subequation $F$ we mean
$$
u+w \ \ {\rm satisfies \ the\ (SMP)\ for\ all\ }\ u\in F(\ob), \ w\in \ft(\ob).
\eqno{(SC)}
$$

Unlike (C), strong comparison (SC) does not always hold for pure second-order constant coefficient
subequations (for instance $F=\cp$). 
In  Section \KK \ we  establish a sufficient condition for (SC) utilizing 
a  ``monotonicity subequation'' $M_F$ associated to $F$.

\Theorem{D} {\sl
If the dual $\wt { M_F}$ satisfies the strong maximum principle, then the
strong comparison principle holds for $F$.
}\medskip

We leave as an open question: When does the strong comparison principle (SC) for
$F$  imply the (SMP) for  $\wt { M_F}$?

Surprisingly, not all monotonicity subequations $M_F$ are convex cones.
In Section \LL,  utilizing such $M_F$,  we construct many new examples of borderline equations for which 
strong comparison holds.  Specifically, for each decreasing continuous function
$g:[0,\infty) \to \bbr$ with $g(0)=0$ and $g(x)<0$ for $x>0$, we construct two
equations $M^g$ and ${\wt M}^g$, with ${\wt M}^g$ borderline, and compute the
characteristic function $f$ of ${\wt M}^g$ in terms of $g$.  

\Theorem{E} {\sl
If $g$ is subadditive and $\int_{0^+} { dy \over f(y)} =\infty$, where $f$ is the characteristic function
associated with ${\wt M}^g$, then the strong comparison principle 
holds for $M^g$ and ${\wt M}^g$.
}
\medskip

Many such functions $g$ exist.  For further examples see (\LL.4) and [\BRUCK], 
 and see Example \LL.7
for a specific example related to the Hopf function (\LL.10).

We use ``increasing'' to mean non-decreasing throughout the paper.

\vfill\eject

%\vskip.3in

%%%%%%%%%%%%%%%%%%%%%%%%%%%%%%%%%%%%%%%%%% 
%%%%%%%%%%%%%%%%%%%%%%%%%%%%%%%%%%%%%%%%%% 
%%%%%%%%%%%%%%%%%%%%%%%%%%%%%%%%%%%%%%%%%% 
%%%%%%%%%%%%%%%%%%%%%%%%%%%%%%%%%%%%%%%%%% 
%%%%%%%%%%%%%%%%%%%%%%%%%%%%%%%%%%%%%%%%%% 

\centerline{\headfont \BB.\  Characterizing the Maximum Principle}.
 \smallskip

In this section we review and amplify the (MP) results in [\DDD].

Let $\cpt$ denote the subset of $A\in \Symn$ with at least one non-negative eigenvalue,
i.e., with $\l_{\rm max}(A)\geq0$.  For the maximum principle we only need to consider
subequations $F\ss \cpt$, since if $A\notin \cpt$, then $A$ is negative definite and $\bra{Ax}x$
violates  (MP).  Note that $\cpt$ is a subequation, that is, it is a closed set which satisfies (P).  
In fact, $\cpt$ is universal for (MP) in the following sense.

\Theorem{\BB.1. (Part I)} \ \  {\sl
Suppose that $F$ is a subequation.
$$
(a)\ \ \ {\rm The \ (MP)\ \ holds\ for\ \ } F 
\qquad\iff\qquad
 F\ss\cpt.  \qquad\qquad\qquad\qquad\qquad\ \ \ \ \ \ \ \ \ 
\eqno{(\BB.1)}
$$
}

\pf It remains to show that (MP) holds for $\cpt$, which follows from Proposition
\BB.3.

\Def{\BB.2} A function $u$ is {\bf subaffine on $X$} if it is upper semi-continuous on $X$ and 
\medskip
\centerline
{
for all compact sets $K\ss X$ and affine functions $a(x) \equiv \bra px +c$,
}
\medskip
\centerline
{
$u\leq a$ \ on\  $\partial K \qquad\Rightarrow \qquad u \leq a $ \ on \ $K$.
}
\bigskip

Subaffine functions clearly satisfy (MP) (take $a(x) = c=$ constant in Definition \BB.2).

Furthermore, for any pure second-order subequation $F$,  {\sl the functions $u\in F(X)$
satisfy (MP) if and only if they are subaffine}, since the sum $u+a$ of a function
$u$ in $F(X)$ and an affine function $a$ is again in $F(X)$.
The subaffine fundtions have an advantage over the larger class of functions satisfying the
(MP) in that they are determined by a local property.

\Prop{\BB.3} 
$$
u \in \cpt(X)\qquad\iff\qquad u\ \ {\sl is\ subaffine\ on\ \ } X.
$$

\pf Suppose $u$ is not subaffine. Then there exists a compact set $K\ss X$ and an affine function
$a$ so that (MP) fails for $w\equiv u-a$ on $K$, i.e.,  $w$ has a strict interior maximum point
on $K$.  This  also holds for $w+\e {|x|^2\over 2}$ with  $\e>0$ sufficiently small.
Then $\vf =-\e {|x|^2\over 2}$  is a test function for $w$ at the maximum point $\bar x
\in\Int K$. Since $D^2_{\bar x}\vf =  - \e I <0$, we conclude that $w\notin \cpt(X)$ 
and so $u\notin\cpt(X)$.

If $u\notin \cpt(X)$, then there exists a test function $\vf$ for $u$ at a point $\bar x\in X$ with
$D^2_{\bar x}\vf \notin \cpt$, i.e., $A\equiv  D^2_{\bar x}\vf <0$.  Set 
$a(x) \equiv \bra {D_{\bar x}\vf}{x-\bar x} +\vf(\bar x)$.  Then   $u(x) \leq a(x) + {1\over 4}\bra {A(x-\bar x)} {x-\bar x}$ near $\bar x$, showing that $u$ is not subaffine on  a small ball $K$ about $\bar x$.\qed

\medskip

In particular, we have, as advertised, that if  $u\in\USC(X)$ is locally subaffine, then $u$ is subaffine.
We will refer to $\cpt$ as the {\sl subaffine subequation}.

Note that in addition to Theorem \BB.1(a)
  we have established the following additional characterizations of the maximum principle.
(The condition in (\BB.3) implies $0\in\Int F$ by  positivity.)
\Theorem {\BB.1. (Continued)}
$$
(b)\ \ \ {\rm The \ (MP)\ \ holds\ for\ \ } F 
\qquad\iff\qquad
0\,\notin\ \Int F.  \qquad\qquad\qquad\qquad\qquad\ \ \ \ \ 
\eqno{(\BB.2)}
$$
$$
(c)\ \ \ {\rm The \  (MP)\ \ fails\ for\ \ } F 
\quad\iff\quad
-\e{|x|^2\over 2} \ \ {\rm is\ \ } F \,{\rm subharmonic \  for\ some\ }\ \e>0.
\eqno{(\BB.3)}
$$

\Remark{\BB.4} Part (a) of Theorem \BB.1 
states that  $\cpt$ is the ``universal''  subequation for  the maximum principle.
Part (b) provides the simplest test for the (MP) to hold for $F$.
Part (c) states that the function  $-\e{|x|^2\over 2}$ is a ``universal'' counterexample to the maximum principle.
\medskip

An obvious corollary of Theorem \BB.1(b) is the following.

\Cor{\BB.5. (Localization)}
$$
\eqalign
{
 & {\sl  If \  two \ subequations\ } F {\sl \ and\ }  G  \ {\sl agree}
  \cr
{\sl    in \ a \ }
& {\sl neighborhood\  of  \  the \  origin \  in \ }  \Symn, \ {\sl then}
\cr
{\rm (MP)} &{\sl \ \ holds\ for\ \ } F 
\qquad \iff\qquad
{\rm (MP)\ \ holds\ for\ \ } G.
}
\eqno{(\BB.4)}
$$
\medskip

A discussion of the subaffine subequation is not complete without 
mentioning its duality with the convex subequation $\cp$.

\vskip.2in
\centerline{\bf Duality}
\bigskip

For any subset $F$ of $\Symn$, the {\bf Dirichlet dual} $\ft$ is defined to be:
$$
\ft\ = \ -\left(\sim\Int F\right) \ =\ \sim\left(-\Int F\right).
\eqno{(\BB.5)}
$$
One can  calculate the key property that
$$
\wt{F+A}\ =\ \ft-A\quad {\rm for \ each\ \ }A\in\Symn.
\eqno{(\BB.6)}
$$
This can be used to show that 
$$
F\ \ {\rm satisfies \ (P)}\quad \Rightarrow\quad  \ft\ \ {\rm satisfies \ (P)}.
\eqno{(\BB.7)}
$$
Other properties of the subequations and their dual subequations include:
$$
F_1\ \ss\ F_2  \ \  \Rightarrow\ \   \ft_2\ \ss\ \ft_1, \qquad
\wt{F_1\cap F_2}\ =\ \ft_1\cup \ft_2, \qquad
\eqno{(\BB.8)}
$$
$$
\wt{\wt F} \ =\ F 
\qquad\qquad
\Int \ft \ =\ -(\sim F)
\qquad\qquad
\partial \ft \ =\ - \partial F
\eqno{(\BB.9)}
$$
The first assertion in (\BB.9) follows from $\Int F \ss \wt{\wt F} \ss F$ combined with 
condition (T) for $F$. The second assertion in (\BB.9) is a restatement of the first.

The Dirichlet dual of $\cp$ can be computed as follows.  Let  $\l_{\rm min}(A)$ and 
$\l_{\rm max}(A)$  denote the smallest and the largest eigenvalues of $A\in\Symn$.
By definition
$$
 \cp\ =\ \{A : \l_{\rm min}(A)\geq0\}.
\eqno{(\BB.10)}
$$
Since $\l_{\rm min}(-A) = -\l_{\rm max}(A)$ it is easy to see that the dual of $\cp$ is
$$
 \cpt\ =\ \{A : \l_{\rm max}(A)\geq0\},
\eqno{(\BB.11)}
$$
justifying the notation $\cpt$ for the subaffine subequation.

\vfill\eject

%%%%%%%%%%%%%%%%%%%%%%%%%%%%%%%%%%%%%%%%%% 
%%%%%%%%%%%%%%%%%%%%%%%%%%%%%%%%%%%%%%%%%% 
%%%%%%%%%%%%%%%%%%%%%%%%%%%%%%%%%%%%%%%%%% 
%%%%%%%%%%%%%%%%%%%%%%%%%%%%%%%%%%%%%%%%%% 
%%%%%%%%%%%%%%%%%%%%%%%%%%%%%%%%%%%%%%%%%% 

\def\Stable{Stable}
\def\stable{stable}

\centerline{\headfont \CC.\  Characterizing the Strong Maximum Principle -- Three Cases}.
 \smallskip

Given a subequation $F$, we consider the following three mutually exclusive cases.

\medskip
{\bf The \Stable \  Case.} \ \ $F \cap (-\cp) \ =\ \emptyset$.

\medskip
{\bf The Algebraic Counterexample Case.} \ \ $(F-\{0\}) \cap (-\cp) \ \neq \ \emptyset$.

\medskip
{\bf The Borderline Case.} \ \ $F \cap (-\cp) \ =\ \{0\}$.

\medskip
\noindent
The first case is stable or generic among subequations where the (SMP) holds,
while in the second case the (SMP) fails via  a quadratic counterexample.
They are both very easy to analyze.

\Theorem{\CC.1}\medskip
{\sl 

(a)\ \ If $F$ is \stable, then the (SMP) holds for $F$ (and for all subequations in a small distance neighborhood of $F$).

\medskip

(b)\ \ If $F$ falls into the algebraic counterexample case, then the (SMP) fails for $F$.
\medskip

(c)\ \ If $F$ is borderline, then the (MP) holds but the (SMP) may or may not hold.
}
\pf
(a)\ \ For completeness first note the equivalent ways of saying that $F$ is \stable.
$$
F\cap (-\cp) \ =\ \emptyset 
\qquad\iff\qquad
F\ \ss\ \Int\cpt 
\qquad\iff\qquad
0\notin F
\eqno{(\CC.1)}
$$
(by positivity, if $A\in F$, $A\leq 0$, then $0\in F$).
Suppose the (SMP) fails for $F$.  Then for some domain $\O$
there exists $u\in F(\ob)$ non-constant, but with an interior maximum point
$x_0$.  The constant function $M\equiv \sup_{\ob} u$ is a test function
for $u$ at $x_0$.  Hence, $0=D_{x_0}^2 u \in F$, so $F$ is not \stable.
The second claim in (a) follows from the last part of (\CC.1).

(b)\ \ By  positivity,
$$
\eqalign
{
(F-\{0\}) \cap (-\cp) \ \neq \ \emptyset
\qquad&\iff\qquad
\exists\, A\leq 0, \ A\neq0 \ \ \ {\rm and}\ \ \  A\in F,\cr
&\iff\quad
-\mu P_e \in F \ \ {\rm for\ some\ \ }  \mu>0\ \ {\rm and}\ \ e\neq 0,
}
\eqno{(\CC.2)}
$$
in which case the functions $u(x) = \half \bra{Ax}x$ and
 $-{\mu\over 2}{\bra ex}^2$ are counterexamples to the (SMP) on any domain
$\O$ containing  the origin.

(c)\ \ The (MP) follows from Theorem \BB.1(b).
Borderline examples where (SMP) holds and where (SMP) fails will be given in Section \JJ \ 
after we prove our main result.\qed
\medskip

 The rest of this section is devoted to further discussion of  the borderline case.

\bigskip
\centerline{\bf Borderline Subequations}
\medskip

There are several equivalent ways of describing  the borderline subequations.

\Lemma{\CC.2} A subequation $F$ is {\bf borderline} if and only if any (or all) of the following 
equivalent conditions holds for $F$.
\medskip
\centerline{(1)  \ \   $0 \in \partial F$\  and\  $F-\{0\} \ss \Int \cpt$.
\qquad\qquad
(1)$'$ \ \   $0 \in \partial   \wt F$\  and\  $\cp-\{0\} \ss \Int \ft$.}

\medskip
\centerline{(2)  \ \   $0 \in \partial F$\  and\  $-\mu P_e \notin F \ \forall\, \mu>0, e\neq0.$ 
\qquad
(2)$'$ \ \   $0 \in \partial   \wt F$\  and\ $\mu P_e \in \Int \wt F \ \forall\, \mu>0, e\neq0.$}

%This condition (\DD.1)$'$ is one definition of $\ft$ being {\bf strictly elliptic at 0}. 
 
\pf 
Since $\Int \cpt = \sim(-\cp)$, (1) is a rephrasing of the definition of borderline.
The equivalences (1) $\iff$ (1)$'$ and 
 (2) $\iff$ (2)$'$ follow from (\BB.8) and (\BB.9). Condition (1) implies Condition (2)
 because $-\mu P_e \notin \Int \cpt$ for $\mu>0$.
Condition (2)$'$ implies Condition (1)$'$  since, by (P), $\Int \wt{F} + \cp \ss \Int \wt{F}$,  and 
$\cp - \{0\}$ is the convex hull of the elements $\mu P_e$ for $\mu>0$ and $e\in \rn$. 
\qed
\medskip

\bigskip
\centerline
{\bf
The Characteristic Function of a Subequation
}
\medskip

In order to further analyze (not necessarily  borderline  or invariant) subequations 
 we  associate two functions $\uf \leq \of$ 
with $F$. We begin by considering a general subequation $F$.
The motivation and more details will be provided later in Section 5.
First we associate the following  two closed sets in $\bbr^2$ with $F$, 
called the {\bf upper (larger) and lower (smaller) radial profiles of $F$}:
$$
\eqalign
{
\oL \ &\equiv \ \left \{(\l,\mu) :  \l P_{e^\perp} + \mu P_{e} \in F \ \ {\rm for\ some\ } e\neq 0\right\},
\cr
\uL \  &\equiv\  \left\{(\l,\mu) :  \l P_{e^\perp} + \mu P_{e} \in F \ \ {\rm for\ all\ } e\neq 0\right\}.
}
\eqno{(\CC.3)}
$$

Since $F$ is $\cp$-monotone,
$$
\oL \ \ {\rm and}\ \  \uL\ \  {\rm are}\ \ \bbr_+\times \bbr_+\ {\rm monotone}.
\eqno{(\CC.4)}
$$

Closed subsets $\L\ss\bbr^2$ which are $\bbr_+\times \bbr_+$-monotone 
can be classified in several ways.
The classification we need is in the following lemma.

\Lemma{\CC.3} {\sl
A set  $\L\ss\bbr^2$ is closed and  $\bbr^2_+$-monotone \ \ \ $\iff$\ \ \ 
there exists a lower semi-continuous,  decreasing function $h:\bbr \to \bbr\cup\{\pm \infty\}$
 such that $\L = \{(\l,\mu) : \mu\geq h(\l)\}$.
}
\pf
Given $\L$, for each $\l\in\bbr$, define $h(\l) = \inf\{\mu : (\l,\mu)\in\L\}$, with
$h(\l)= - \infty$ if this set is all of $\bbr$ and $h(\l)=\infty$ if this set is empty.
The $\bbr^2_+$-monotonicity implies that $h$ is decreasing.  Now $\L$ is closed if and only
if $h$ is lower semi-continuous. The remainder of the proof is left to the reader.\qed
\medskip

It is more convenient to replace $h$ by the function $f\equiv-h$ so that 
$f:\bbr \to \bbr\cup \{\pm\infty\}$ is upper semi-continuous, increasing and 
$$
\L \ \equiv\ \{(\l, \mu) : \mu+f(\l)\geq 0\}.
\eqno{(\CC.5)}
$$
Thus the radial profiles $\oL$ and $\uL$  of $F$ can be used interchangeably with the 
following associated
functions $\of$ and $\uf$ describing them.

\Def {\CC.4}  
Suppose that $F$ is a   subequation.
The {\bf upper (larger)  and lower (smaller)  characteristic functions $\of$ and $\uf$
associated with $F$} are defined  by:
$$
\of(\l) \ \equiv \ \sup\left\{ \mu : \l P_{e^\perp} - \mu P_{e} \in F \ \ {\rm for\ some\ } e\neq 0\right\}
$$
$$
\uf(\l) \ \equiv \ \sup\left\{ \mu : \l P_{e^\perp} - \mu P_{e} \in F \ \ {\rm for\ all\ } e\neq 0\right\}
$$

Summarizing, we have 
$$
\l P_{e^\perp} + \mu P_{e} \in F \ \ \ {\rm for\ some\ \ }    e\neq 0
\qquad\iff\qquad \mu+ \of(\l)\geq 0.
\eqno{(\CC.6)}
$$
$$
\l P_{e^\perp} + \mu P_{e} \in F \ \ \ {\rm for\ all\ \ }    e\neq 0
\qquad\iff\qquad \mu+ \uf(\l)\geq 0.
\eqno{(\CC.7)}
$$

We will use the following fact to further analyze the borderline case.

\Lemma {\CC.5} 

$$
F\ \ {\rm is\  borderline} \qquad\iff\qquad \uf(0) \ =\ \of(0)\ =\ 0.
\eqno{(\CC.8)}
$$
 
 \pf Use Definition \CC.4 and  condition (2)  in Lemma \CC.2.\qed
\medskip

 The asymptotic structure of $F$ near 0 is reflected 
in the asymptotic behavior of $\uf$ and $\of$ near 0.
Now we can state the main result of this paper.
Note that only the behavior of $\uf(\l)$ and $\of(\l)$ for $\l$ positive (and small) affects the outcomes.

\Theorem{\CC.6} {\sl Suppose that $F$ is a borderline subequation with upper and lower characteristic
functions $\of$ and $\uf$.
$$
{\rm (a)} \quad {\sl If}\ \ \int_{0^+} {dy\over \of(y)}\ =\   \infty, \ \ {\sl then\  the \ (SMP)\ holds\ for\ }F. 
$$
$$
{\rm (b)} \quad {\sl If}\ \ \int_{0^+} {dy\over \uf(y)}\ <\  \infty, \ \ {\sl then\  the \ (SMP)\ fails\ for\ } F.
 $$
}
\medskip

The only case not covered is when $ \int_{0^+} {dy\over \of(y)} < \infty$
and $\int_{0^+} {dy\over \uf(y)}\ =  \infty$.
In this case  $F$  will be referred to as  a {\bf gap subequation}.

\Def{\CC.7. (Weakly Invariant Subequations)}  
For most equations of interest, $\uf = \of$, and in this case Theorem \CC.6
 gives a necessary and sufficient condition
for $F$ to satisfy the (SMP). First 
note that $\uf =\of$ if and only if $\uL=\oL$, or equivalently, for all $\l,\mu$:
$$
{\rm If\ } \l P_{e^\perp} +\mu P_e \in F \ \ {\rm for\ some\ \ } e\ne0, \ \ {\rm then\ } 
 \l P_{e^\perp} +\mu P_e \in F \ \  {\rm for\ all\ \ } e\ne0
\eqno{(\CC.9)}
$$
We take this as the definition of $F$ being {\bf weakly invariant}, and for simplicity
we shall refer to it by just saying  that $F$ is {\bf invariant}.
Note also that $P_{e^\perp}, P_e$ have the same span as $I,P_e$, and therefore, 
for any subequation $F$ which
 is invariant under the action of a group $G$ acting transitively on the
unit sphere $S^{n-1}\ss\rn$, the characteristic functions $\uf$ and $\of$ are equal.
  Among possibilities for  $G$ 
are SO$_n$ acting on $\rn$, SU$_n$ acting on $\bbr^{2n}=\bbc^n$, 
Sp$_n$ acting on $\bbr^{4n}=\bbh^n$, G$_2$ acting on $\bbr^7$ and Spin$_7$ acting
on $\bbr^8$.

\def\psii{\chi}
\medskip
\noindent
{\bf  Some Examples \CC.8.}  Suppose $\psii :\bbr \to \bbr$ is odd ($\psii(-t)=-\psii(t)$) 
and strictly increasing. Fix $1\leq k\leq n$ and define  $F=F_{\psii, k}$ to
be the set of $A\in\Symn$ such that 
$$
\s_\ell(\psii(A)) \ \geq\ 0\qquad {\rm for}\  \ell=1,...,k,
$$
where $\s_\ell$ denotes the $\ell^{\rm th}$ elementary symmetric function.  That is,
$A\in F$ if and only if 
$$
\s_\ell\left( \psii(\l_1(A)),...,\psii(\l_n(A))    \right) \ \geq\ 0 \qquad {\rm for}\  \ell=1,...,k,
$$
where $\l_1(A),...,\l_n(A)$ are the eigenvalues of $A$. One checks that $F$ satisfies 
Condition (P) and is therefore a subequation.  Direct calculation shows that the characteristic 
function $f=\of=\uf$ is
$$
f(\l)\ =\ \psii^{-1}\left\{\left({n\over k} - 1\right) \psii(\l)\right\}.
\eqno{(\CC.10)}
$$

One concludes that{\sl  if $\psii$ is smooth, or just Lipschitz,  in a neighborhood of 0, then the (SMP) holds.}
\medskip

A basic case is where  $\psii(t) = {\rm sign}(t) |t|^\b$ for $\b>0$.
For example, if  $k=1$ and $\b={1\over 3}$,  then $F=\{A: \tr(A^{1\over 3})\geq0\}$.
In these cases $f(\l) = ({n\over k}-1)^{1\over \b} \l$, and so the (SMP) holds.

In fact, by Theorem \JJ.5 below, this basic example is contained in the cone
subequation $\cp_{\a}^{\rm min/max}$ with $\a \equiv ({n\over k}-1)^{1\over \b}$.
Now subequations which are cones can  be treated more classically.
Nevertheless, they are useful in understanding our main result in the non-invariant case,
so we examine them next.

\bigskip
\centerline{
\bf Local Cone Subequations}
\medskip

Perhaps the simplest examples where Theorem \CC.6 applies are the cone subequations.
The results in this case are  not really new, but they give a nice illustration of our geometric
point of view   and parallel the characterizations  obtained for the (MP) as described in Remark \BB.4.

We say that a subset $F\ss\Symn$  is a {\bf cone} if $tF\ss F$ for all $t>0$, and a {\bf local cone}
if for some $\d>0$ we have that $F\cap B_\d(0)$ is the cone on the (non-empty) link $F\cap \partial B_\d(0)$. 

\Theorem{\CC.9. ($n\geq2$)}
{\sl
Suppose that $F$ is a local cone subequation.
\medskip
\noindent
(a) \ \ The (SMP) holds for $F \quad \iff \quad F-\{0\}\ss \Int\cpt.$
\medskip
\noindent
(b) \ \ The (SMP) holds for $F \quad \iff \quad - \e P_e\notin F \ \ \forall\, \e>0 \ \ {\rm and} \ \ |e|=1.$
\medskip
\noindent
(c) \ \ The (SMP) fails for $F \quad \iff \quad -\smfrac\e 2{\bra e x}^2$  is $F$-subharmonic for some 
$\e > 0$ and 

\hskip 3in and some $|e|=1$.

\medskip
Moreover, 
$$
\int_{0^+} {dy\over \of(y)}= \infty \ \ {\sl
  is\  both\  necessary\  and\  sufficient \ for\  the \ (SMP)}
\eqno{(\CC.11)}
$$
(cf.  Theorem  \CC.6(a)),
%Moreover, the condition $\int_{0^+}(1/\of) = \infty$ 
%is both necessary and sufficient for the (SMP)
and is equivalent to $F$ being borderline.
}
\pf
For a local cone subequation we have $0\in\partial F$. Hence by Lemma \CC.2, part (2)
$$
F \ \ {\rm is\ borderline}\ \ \quad\iff\quad - \e P_e \notin F \ \ \forall\, \e>0 \ \ {\rm and}\ \ e\neq 0.
\eqno{(\CC.12)}
$$
Said differently, $F$ is not borderline $\iff - \e P_e \in F$ for some $\e>0$ and $e\neq 0$.  This proves that 
$$
{\rm either}\ \ F \  {\rm is\ borderline}\ \ {\rm or}\ \  h(x) \equiv - \smfrac{\e}2 { \bra ex}^2
=  - \smfrac{\e}2 \bra {P_ex} x \ \ {\rm is}\ F\, {\rm subharmonic},
\eqno{(\CC.13)}
$$
in which case the (SMP) fails for $F$.

Define $0\leq\overline \a \leq \infty$ by
$$
\overline \a \ \equiv\ \sup  \left\{\a : % {\d \over \sqrt{n-1+\a^2}}
 \smfrac {\d}{\sqrt{n-1+\a^2}}( P_{e^\perp} - \a P_e) \in \partial B_\d(0)\cap F\ \ {\rm for\ some\ } 
|e|=1  \right\}
\eqno{(\CC.14)}
$$
The local cone condition  implies that if $\overline \a<\infty$, then
$$
\of (\l) \ =\ \overline{\a} \l \ \ \ {\rm for\ \ } \l \ \leq\  \l_0
\eqno{(\CC.15)}
$$
 where $\l_0 = ({1\over 1+\a^2})^{1\over2}$.
Hence, $\int_{0^+} {1\over \of} = {1\over \overline\a} \int_{0^+} {d\l\over \l} =\infty$,
 in which case by Theorem \CC.6(a), $F$ satisfies the (SMP).
On the other hand, if
 $\overline \a=\infty$, then since $F\cap \partial B_\d(0)$ is closed, 
 one has that $-\d P_e \in F$ for some  $|e|=1$, and hence the (SMP) fails.

Finally, note that
conditions in (a) and (b) are just two ways of saying that $F$ is borderline (cf.  Theorem  \CC.6(a)).
In addition, this proves (\CC.11). \qed
\medskip

The Remark \BB.4 describing Parts (a), (b) and (c) of Theorem \BB.1 has the following parallel 
describing parts (a), (b) and (c) of Theorem \CC.9.

\Remark{\CC.10. (The (SMP) for Local Cones)}
 Part (a) states that $\Int \cpt$ is a ``universal'' set for the (SMP), while part (b)
gives the simplest test for the (SMP).
Part (c) says that one-variable quadratic functions such as $-\half x_1^2$ provide a 
``universal'' set of counterexamples to the (SMP).\smallskip

 \Remark{\CC.11} 
The hard half of Theorem \CC.9 (a) or (b) says that borderline cone subequations satisfy the (SMP).
This is easy to prove using the classical ``Hopf Lemma'' construction, and so, in this sense,
Theorem \CC.9 is not new.

\Ex{\CC.12. (Gap Cone Subequations)}  
These are cone subequations where $\int_{0^+}(1/\of) = 0$ (equivalently, $\overline \a =\infty$)
and $\int_{0^+}(1/\uf) = \infty$ (equivalently, $\underline \a <\infty$).
Neither Part (a) nor Part (b) of Theorem \CC.6 applies.  However, Theorem \CC.9 does apply (see (\CC.11))
to show that the (SMP) fails for all gap cone subequations.

Such subequations are easy to construct.  For example, on $\bbr^2$ define $F$ by:
$u_{xx} \geq0$ or, if $u_{xx} <0$, then $u_{yy} + \underline \a u_{xx} \geq0$.

\vfill\eject
%\vskip .3in

%%%%%%%%%%%%%%%%%%%%%%%%%%%%%%%%%%%%%%%%%% 
%%%%%%%%%%%%%%%%%%%%%%%%%%%%%%%%%%%%%%%%%% 
%%%%%%%%%%%%%%%%%%%%%%%%%%%%%%%%%%%%%%%%%% 
%%%%%%%%%%%%%%%%%%%%%%%%%%%%%%%%%%%%%%%%%% 
%%%%%%%%%%%%%%%%%%%%%%%%%%%%%%%%%%%%%%%%%% 

\centerline{\headfont \EE.\ The Radial Subequation Associated to $F$.}
 \smallskip 

Supppose $\psi$ is of class $C^2$ on an interval contained in the positive real numbers.
Consider $\psi(|x|)$ as a function on the corresponding annular region in $\rn$.

\Lemma{\EE.1} {\sl
$$
D^2_x \psi \ =\ {\psi'(|x|)\over |x|} P_{x^\perp} + \psi''(|x|) P_{x}.
$$
%where $P_{[x]} = {x\circ x \over |x|^2}$ 
%denotes orthogonal projection onto the line $[x]$ 
%through $x\neq0$ and $P_{[x]^\perp} =I- P_{[x]}$ denotes orthogonal projection onto
%the hyperplane with normal $[x]$.
}
\pf
First note that $D(|x|)= {x\over |x|}$, and therefore  $D^2(|x|) = D({x\over |x|})
 = {1\over |x|} I - {x\over |x|^2}\circ {x\over |x|}
= {1\over |x|}(I- P_{x} ) = {1\over |x|} P_{x^\perp}$.
Hence, 
$$
D_x\psi \  =\  \psi'(|x|) {x\over |x|} 
\qquad {\rm and}
$$ 
$$D_x^2\psi \ \ =\  \ \psi'(|x|) 
D\left( {x\over |x|} \right) + \psi''(|x|)  {x\over |x|}\circ  {x\over |x|}
\ \ =\  \ 
{\psi'(|x|) \over |x|} P_{x^\perp} + \psi''(|x|)P_{x}. \qquad\qquad
\mathqed
$$

\Cor{\EE.2}  {\sl
The second derivative $D^2_x \psi$ has  eigenvalues ${\psi'(|x|) \over |x|}$
with multiplicity $n-1$ and   $ \psi''(|x|)$ 
with multiplicity 1.
}\medskip

For simplicity we shall now assume that $F$ is invariant as in 
Definition \CC.7 and let $f=\of=\uf$ denote its characteristic function.
  Recall from (\CC.6), or (\CC.7),  that
 $$
 \l  P_{e^\perp} + \mu  P_{e} \ \in\ F \ \ \forall \, e\neq0
 \qquad\iff\qquad
 \mu+  f(\l)\geq0.
\eqno{(\EE.1)}
$$

With motivation from Lemma \EE.1 this leads to a subequation $R_f$ on $(0, \infty)$.
Let $p=\psi'(t)$ and $a=\psi''(t)$ denote jet coordinates.

\Def{\EE.3}  The {\bf radial subequation $R_f$ associated to $F$} is defined by
$$
R_f \ :\  a+f \left({p\over t}\right)\ \geq\ 0\qquad 0<t<\infty
\eqno{(\EE.2)}
$$
where $f$ is the characteristic function associated with the subequation $F$.

\medskip

It follows immediately from these definitions and Lemma \EE.1 that if $\psi(t)$ is a 
$C^2$-function defined on a subinterval of $(0,\infty)$,  with $u(x)\equiv \psi(|x|)$ defined on the corresponding 
annular region in $\rn$, then
$$
u(x) \equiv \psi(|x|) \ \ {\rm is}\ \ F\ {\rm subharmonic}
\quad   \iff    \quad
\psi(t) \ \ {\rm is}\ R_f \ {\rm subharmonic} 
\eqno{(\EE.3)}
$$
  
This is extended  to upper semi-continuous  functions in  Appendix  A (Theorem A.1). 
The proof of the implication $\Rightarrow$ is elementary, whereas the proof of
$\Leftarrow$  requires some details. However, note that $u(x) = \psi(|x|)$ is upper semicontinuous $\iff$
$\psi(t)$  is upper semicontinuous.

\Remark{\EE.4}
The radial subequation $R_f$ associated to $F$ satisfies the topological conditions
(T) in [\DDR].  Namely, 
$$
(i)\ \  R\ =\ \overline{\Int R},  \qquad
(ii)  \ \ R_t\ =\ \overline{\Int R_t},  \qquad
(iii)  \ \ \Int_t R_t\ =\  (\Int R)_t.
$$
where $R_t$ is the fibre of $R$ above $t$ and $\Int_t$ denotes the interior in $R_t$.
Note that $\Int R$ is not defined by $a+f({p\over t})>0$ but by 
$a+f_-({p\over t})>0$ where $f_-(y) \equiv \lim_{z\to y^- }f(z)$ is lower semi-continuous.
The proof is left to the reader.

\Remark{\EE.5. (Radial Harmonics)}
If $\psi(t)$ is $R_f$-harmonic on an interval $I\ss (0,\infty)$, then for any constants
$r>0$ and $k\in \bbr$, the 2-parameter family of functions
 $\psi_r(t) \equiv r^2\psi(t/r) + k$ consists of  $R_f$-harmonics on $rI$.
This follows since $\vf$ is a test function for $\pm \psi$ if and only if $\vf_r$ is a test function 
for $\pm \psi_r$, and the assertion is true for $C^2$-functions.

 %\vfill\eject
\vskip.3in

%%%%%%%%%%%%%%%%%%%%%%%%%%%%%%%%%%%%%%%%%% 
%%%%%%%%%%%%%%%%%%%%%%%%%%%%%%%%%%%%%%%%%% 
%%%%%%%%%%%%%%%%%%%%%%%%%%%%%%%%%%%%%%%%%% 
%%%%%%%%%%%%%%%%%%%%%%%%%%%%%%%%%%%%%%%%%% 
%%%%%%%%%%%%%%%%%%%%%%%%%%%%%%%%%%%%%%%%%% 

\centerline{\headfont \GG.\  Increasing Radial Subharmonics for Borderline Subequations}.
 \smallskip
 
 As in the last section, we assume for simplicity that $F$ is invariant, i.e., $f=\uf=\of$.  
 Because of the next result we  focus on radial subharmonics which are increasing.
  
 \Lemma {\GG.1} {\sl
 Suppose that $F$ is borderline  and $u(x) =\psi(|x|)$ is a radial $F$-subharmonic function.
 There is only one way that $u(x)$ can violate the (SMP).  Namely, for some $r$,
 $\psi(t)$ must satisfy:
 $$
 \psi(t)\ <\ M\ \ {\sl for}\ \ t<r
 \and
  \psi(t)\ \equiv\ M\ \ {\sl for}\ \ t\geq r.
 \eqno{(\GG.1)}
 $$
 Moreover,}  
  $$
\psi \ {\sl must \  be \  increasing\ on\ } (\bar a, r), \  {\sl for \  some\ }  \bar a < r.
  \eqno{(\GG.2)}
 $$

 \pf
By the borderline hypothesis the (MP) holds for $F$.
 Since $u$ satisfies the (MP), so does $\psi(t)$. If $\psi$ has an interior maximum
 point at $t_0$ on an interval $[a,b]$, then either $\psi$ is equal to the maximum value
 $M$  on $[a,t_0]$
 or on $[t_0,b]$ since otherwise $\psi$ violates the (MP) on an interval about $t_0$.
 If $\psi$ equals this maximum value on $[a,t_0]$, we can extend $u(x)$ to the ball of 
 radius $b$ (to be constant  on the ball of radius $a$) as an $F$-subharmonic function
 which violates the (MP).  This proves (\GG.1).
 
 Pick $\bar a$ to be a minimum point for $\psi$ on $[a,b]$.
 Then $\psi$ must be increasing on $[\bar a, b]$. Otherwise, there exist $\a,\b$ with
 $\bar a < \a < \b < r$ and $\psi(\a) > \psi(\b)$.  In this case $\psi(\a) > \psi(\bar a)$ also 
(since $\psi(\b)\geq\psi(\bar a)$),  and this violates  the maximum principle on  $[\bar a, \b]$.
\qed

\def\ua{\uparrow}
\def\da{\downarrow}

\Def{\GG.2}  Suppose that $F$ is borderline.   
The {\bf  increasing radial subharmonic equation}  $R_f^\ua$ on $(0,\infty)$ is
defined by
$$
R_f^\ua   \  :\   a+f\left( {p\over t}\right) \ \geq\ 0 \ \ {\rm and}\ \ p\geq0. 
\eqno{(\GG.3)}
$$
where $f$ is the characteristic function of $F$.\medskip

For $C^2$-functions $\psi(t)$ it is obvious that:
$$
\psi(t) \ \ {\rm is \ } R_f^\ua \, {\rm subharmonic}  \quad \iff\quad \psi(|x|) \ \ {\rm is \ } F\cap\{x\cdot p\geq0\}  \,{\rm subharmonic}
\eqno{(\GG.4)}
$$
where $F\cap \{x \cdot p\geq0\}$ is a variable coefficient subequation on $\rn$
depending on both the first and second derivatives.
The equivalence  (\GG.4) is  extended using Theorem A.1.

\Theorem{\GG.3. (Increasing Radial Subharmonics)}  {\sl
Suppose that $F$ is borderline.  The function $u(x) =\psi(|x|)$ is
$F$-subharmonic and radially increasing on an annular region in $\rn$
if and only if  $\psi(t)$  is $R_f^\ua$-subharmonic on the corresponding subinterval of $(0,\infty)$.}

\pf  Theorem A.1 states that: $u$ is $F$-subharmonic $\iff$ $\psi$ is $R_f$-subharmonic.
By definition $u(x)$ is radially increasing if $u$ satisfies the first-order variable coefficient subequation
$\{p\cdot x\geq0\}$. It remains to show that
$$
u\ \ {\rm satisfies\ the\ subequation} \ \{p\cdot x\geq0\}
\qquad\iff\qquad 
\psi\ \ {\rm satisfies
 } \ \psi'(t)\geq0.
\eqno{(\GG.5)}
$$

 Suppose  $\psi(|x|)$ is  $\{x\cdot p\geq0\}$-subharmonic and that $\vf(t)$ is a 
test function for $\psi(t)$ at a point $t_0$.  Then $\vf(|x|)$ is a test function for $u(x)$  at $x_0$
if $|x_0|=t_0$.  Now
$$
D_{x_0}  \vf \ =\ \vf'(|x_0|) {x_0\over |x_0|}\quad {\rm and\ hence}\quad
x_0\cdot D_{x_0}\vf \ =\ |x_0|\vf'(|x_0|).
\eqno{(\GG.6)}
$$
Thus $\vf'(t_0)\geq 0$ proving that $\psi(t)$ is increasing.  Conversely,
if $\psi(t)$ is increasing and $\vf(x)$ is a test function for $u(x)$ at $x_0$,
then $\overline \vf (t) \equiv \vf({tx_0\over |x_0|})$ is a test function for $\psi(t)$ at $t_0=|x_0|$.
Hence, ${\overline \vf }'(t_0)\geq 0$.  However, ${\overline \vf }'(t_0)= (D_{x_0}\vf)\cdot x_0$.\qed

\Remark{\GG.4. (Decreasing Radial Subharmonics)}
For borderline $F$ we define the {\bf decreasing radial subharmonic equation} 
 $R_f^{\da}$ on $(0,\infty)$ by
$$
R_f^{\da}   \  :\   a+f\left( {p\over t}\right) \ \geq\ 0 \ \ {\rm and}\ \ p\leq0. 
\eqno{(\GG.7)}
$$
where again $f$ is the characteristic function of $F$.
We leave it to the reader to show the following.  For  $\psi$ upper semi-continuous,
$$
\psi(t) \ \ {\rm is \ } R_f^\da \, {\rm subharmonic}  \quad \iff\quad \psi(|x|) \ \ {\rm is \ } F\cap\{x\cdot p\leq0\}  \,{\rm subharmonic}
\eqno{(\GG.8)}
$$

% \vfill\eject
\vskip.3in

%%%%%%%%%%%%%%%%%%%%%%%%%%%%%%%%%%%%%%%%%% 
%%%%%%%%%%%%%%%%%%%%%%%%%%%%%%%%%%%%%%%%%% 
%%%%%%%%%%%%%%%%%%%%%%%%%%%%%%%%%%%%%%%%%% 
%%%%%%%%%%%%%%%%%%%%%%%%%%%%%%%%%%%%%%%%%% 
%%%%%%%%%%%%%%%%%%%%%%%%%%%%%%%%%%%%%%%%%% 

\centerline{\headfont \HH.\ Proof of the (SMP)}.
 \medskip

\def\Fh{F^{\#}}

 In this section we prove Part (a) of Theorem \CC.6.
 The subequation $F$ is assumed to be borderline, and 
 we can assume that it  is O$_n$-invariant because 
 of the following construction.  Set
 $$
 \Fh \ \equiv\ \bigcup_{g\in {\rm O}_n} g(F).
 \eqno{(\HH.1)}
 $$
 First note that $ \Fh$ is also a subequation.
Now from Definition \CC.4 of the characteristic function $\of$ of $F$ and the fact
that $P_{e^\perp}$, $P_e$ have the same span as $I$, $P_e$, it is easy to see that
the characteristic function for $ \Fh$ is $\of$.  Moreover, $ \Fh$ is an O$_n$-invariant 
subequation which contains $F$ so that it suffices to prove Theorem \CC.6(a) for $ \Fh$.

From now on we assume that  $F$ is an  O$_n$-invariant borderline subequation,
and we let $f$ denote the restriction of $\of=\uf$ to $[0,\infty)$.  Hence $f(0)=0$
and $f$ is increasing.  Furthermore, let both $R_f^\ua$ and  $R_F^\ua$ denote the
subequation defined by (\FF.4).

Part (a) of Theorem \CC.6 follows from two implications.

\Lemma {\HH.1}
{\sl Suppose $F$ is an  O$_n$-invariant borderline subequation.  Then}
$$
\int_{0^+} {dy\over f(y)}\ =\ \infty 
\qquad\Rightarrow\qquad
{\rm (SMP)\ \ for\ \ } R^\ua_f,\ \ {\rm and}
\eqno{(\HH.2)}
$$
$$
{\rm (SMP)\ \ for\ \ } R^\ua_F
\qquad\Rightarrow\qquad
{\rm (SMP)\ \ for\ \ } F
\qquad\ \ \
\eqno{(\HH.3)}
$$
\pf
We prove the second implication (\HH.3) first. 
Suppose $u$ is a counterexample  to the (SMP) for $F$ on a bounded domain $\O$.
We will show this
leads to a counterexample to the (SMP) for $R_F^\ua$.

First, for all sufficiently small $\overline r>0$,  there exists a ball $B_{\overline r}(x_0) \ss\O$ of radius
$\overline r$ such that the maximum $M\equiv \sup_{\ob} u$ satisfies
$$
\eqalign
{
&(a)\ \ \ u(x) \ <\ M\quad {\rm for\ all\ }\ \ x\in B_{\overline r}(x_0) 
\and  \cr
&(b)\ \ \ u(\overline x) \ =\ M\quad {\rm for\ some\ }\ \ \overline x\in \partial B_{\overline r}(x_0).
}
\eqno{(\HH.4)}
$$
This can be seen as follows.  Since (SMP) is false, there exist points in $\O$ which are not
in the maximum set $\{u=M\}$.  Pick such a point $x_0$ closer to $\{u=M\}$ than to $\bo$
and set $\overline r \equiv {\rm dist}(x_0, \{u=M\})$.  Let $B_t \equiv B_t(x_0)$ and
$M(t) \equiv \sup_{\partial B_t} u$.

Second, choose an annulus
$$
A\ =\ A(r,R)\ \equiv\ \{x : r \ \leq\ |x-x_0|\ \leq\ R\} \ \ss \ \O\ 
\eqno{(\HH.5)}
$$
containing $\partial B_{\overline r}$ in its interior. i.e., with $r<\overline r < R$.
Then
$$
u(\overline x)\ =\ M \ \ {\rm at}\ \ \overline x \in \Int A,\ \ {\rm while \ on\ }\partial A: \ \ u\bigr|_{\partial B_r} \ <\ M\ \ {\rm and}\ \ u\bigr|_{\partial B_R} \ \leq\ M.
\eqno{(\HH.6)}
$$

Since $F$ is borderline, $0\in \partial F$, and hence by Theorem  \BB.1 the (MP) holds for 
$u$ on $B_t$ since $u$ is $F$-subharmonic on $\O$.  Therefore $M(t)$ must be increasing for
$r <t< R$.  Hence, by (\HH.4) and (\HH.6)
$$
M(t)\ < \ M  \quad {\rm for\ \ } \ \  r <t<  \bar r  \and   M(t)=M
\ \  {\rm for\  \  }   \bar r \leq t<   R.
\eqno{(\HH.7)}
$$
That is, the (SMP) for $M(t)$ on $r  \leq t\leq R$ fails.
It remains to show that $M(t)$ is $R_F^\ua$-subharmonic.

\Lemma {\HH.2} {\sl  
For any upper semi-continuous function $u$, the function $M(t) \equiv \sup_{\partial B_t} u$ is upper
semi-continuous.
}
\pf
Assume the balls $B_t$ are centered at the origin.  Given $\d>0$,
$$
N_\d \ \equiv\ \{x : u(x) < M(t)+\d\}
$$
is an open set containing $\partial B_t = \{x:|x|=t\}$.  Hence the annulus $\{x: t-\e\leq |x|\leq t+\e\}$
is contained in $N_\d$ for $\e>0$ small.  Thus $M(r) < M(t) +\d$ if  $t-\e\leq  r\leq t+\e$.
This proves that $M(t)$ is upper semi-continuous.\qed
\medskip

Since $M(t)$ satisfies the subequation $\{p\geq0\}$ it remains to show that $M(t)$
satisfies the subequation $R_F$.  By Theorem A.1  it suffices to show that 
$M(|x|)$ is $F$-subharmonic on $r <|x|<  R$.  The next result completes the 
proof of  (\HH.3).

\Lemma{\HH.3}  {\sl
If $u$ is $F$-subharmonic on an annulus, then $M(|x|)$ is also $F$-subharmonic 
on the same annulus where $M(t) \equiv \sup_{|x|=t} u$.
}

\pf  By Lemma \HH.2 $M(t)$ is upper semi-continuous, and hence $M(|x|)$ is 
upper semi-continuous.  Let $u_g(x)  \equiv  u(gx)$ with $g\in  {\rm O}_n$.
Each $u_g$ is $F$-subharmonic since $F$ is O$_n$-invariant.  Thus
$$
M(|x|)\ =\ \sup_{g\in {\rm O}_n}  u_g(x)
\eqno{(\HH.6)}
$$
is  $F$-subharmonic by the standard ``families locally bounded above'' property for $F$.
\qed

\vskip.3in

%%%%%%%%%%%%%%%%%%%%%%%%%%%%%%%%%%%%%%%%%% 
%%%%%%%%%%%%%%%%%%%%%%%%%%%%%%%%%%%%%%%%%% 
%%%%%%%%%%%%%%%%%%%%%%%%%%%%%%%%%%%%%%%%%% 
%%%%%%%%%%%%%%%%%%%%%%%%%%%%%%%%%%%%%%%%%% 
%%%%%%%%%%%%%%%%%%%%%%%%%%%%%%%%%%%%%%%%%% 

\centerline{\bf A One-Variable Result}.
 \smallskip

The point of this subsection is to prove the one-variable result (\HH.2) 
which completes the proof of Theorem \CC.6 part (a).
We assume throughout that $f:[0,\infty)\to [0,\infty]$ is an upper semi-continuous, increasing
function  with $f(0)=0$, and we define the subequation
$R_f^\ua$ on $(0,\infty)$ by (\GG.4).

\Prop{\HH.4}
$$
\int_{0^+} {dy\over f(y)} \ =\ \infty\qquad \Rightarrow\qquad {\sl The \ (SMP)\ holds\ for\ } R_f^\ua.
$$
\medskip

To prove this we first consider the following one-variable constant coefficient subequation $E$
defined by 
$$
E \ : \ \ a+f(p)\ \geq\ 0 \ \ \ {\rm and}\ \ \ p\geq0.
\eqno{(\HH.7)}
$$

\Prop{\HH.5} {\sl
$$
\int_{0^+} {dy\over f(y)} \ =\ \infty\qquad \Rightarrow\qquad {\sl The \ (SMP)\ holds\ for\ } E.
$$
}

\medskip
\noindent
{\bf Proof that \ \HH.5 \ $\Rightarrow$\ \HH.4.}
  Suppose that the (SMP) fails for $R_f^\ua$ on $[r_1, r_2] \ss (0,\infty)$. 
Choose $r$ with $0<r < r_1$.
Consider the constant coefficient subequation $E_r$ defined by
$$
E_r  \ :\ a+f \left( {p\over r }\right ) \ \geq\ 0\and p\geq0.
\eqno{(\HH.8)}
$$
If $t>r$, then $ a+f \left( {p\over t}\right ) \geq0$ implies that $ a+f \left( {p\over r }\right ) \geq 0$
since $f$ is increasing.  That is, each fibre $(R_f^\ua)_t \ss E_r$ if $t>r$,  so that
on a neighborhood of $[r_1, r_2]$, if $\psi$ is $R_f^\ua$-subharmonic, then $\psi$ is $E_r$-subharmonic.
 Therefore the (SMP) fails for $E_r$.
The function $f({y\over r})$ satisfies the same conditions as the function $f$.
Hence, by Proposition \HH.5,
$
\int_{0^+} {dy\over f(y)} \ =\ {1\over r} \int_{0^+} {dy\over f({y\over r})} \ <\ \infty.
$
\qed

\medskip
\noindent
{\bf Proof of Proposition \HH.5.}
Suppose that $\psi$ is a counterexample to the  (SMP) for $E$.
Since $\psi$ is upper semi-continuous and increasing, there exists a point
$r_0$ such that 
$$
\psi(t) \ <\  M\ \ \ {\rm for}\ \ t<r_0, \qquad {\rm and }\qquad \psi(t) \ \equiv\ M
\ \ \ {\rm for}\ \ r_0 \leq t.
\eqno{(\HH.9)}
$$
By sup-convolution we may assume that $\psi$ is quasi-convex and still satisfies
$E$ with a new $r_0$ slightly smaller than the old one.
Since $f$ is increasing we have the following.

\Lemma{\HH.6} {\sl
The derivative $\psi'$ can be assumed to be absolutely continuous.
}

\pf
Since $\psi(t)+\half \l t^2$ is convex for some $\l>0$, the second distributional
derivative $\psi''=\mu-\l$ where $\mu \geq 0$ is a non-negative measure.
Consider the Lebesgue decomposition $\mu=\a+\nu$ of $\mu$ 
into its $\lloc$-part $\a$ and its singular part $\nu$.
Since $\nu$ is supported on $t\leq r_0$, there exists a unique convex function
$\b$ with $\b''=\nu$ and $\b\equiv 0$ on $r_0\leq t$.   It follows easily that 
$\b(t) \geq0$ and $\b$ is decreasing.  Therefore
$\bar{\psi} (t) \equiv \psi(t) -\b(t)\leq \psi(t)$ and $\bar{\psi} (t)$ is increasing.
Hence $\bar{\psi}$ also satisfies (\HH.9).  Now $\bar{\psi}'' = \a-\l$,
and therefore $\bar{\psi}'$ is absolutely continuous.  Since $\nu$  is 
singular, $\b''(t) = 0$  a.e., and since $\b$ is decreasing, 
$\bar{\psi}'(t) = \psi'(t) - \b'(t) \geq \psi'(t)$ a.e..
Therefore, since $f$ is increasing and $\psi$ is $E$-subharmonic,
$$
\bar{\psi}''(t) + f(\bar{\psi}'(t))\ \geq\ 0\qquad{\rm a.e.}
\eqno{(\HH.10)}
$$
This almost-everywhere inequality is all that will be used to complete the 
proof of Proposition \HH.5.  However, in general, if a quasi-convex function
satisfies a subequation $F$ a.e., then it must be $F$-subharmonic (see Corollary 7.5 in [\DDD]
for pure second-order case and (\II.3) below for the general case).\qed
\medskip

Now  let $\vf(t) \equiv \psi'(t)$.  
This function $\vf$ is absolutely continuous since $\vf'(t) \equiv \a(t)-\l$.
The properties that $\psi$ is increasing and $\psi(t) \equiv M$ for $t\geq r_0$
translate into the properties:
$$
\vf(t) \ \geq\ 0\and \vf(t) \ =\ 0 \ \ {\rm if\ \ } t\geq r_0.
\eqno{(\HH.11)}
$$
The inequality (\HH.10) states that
$$
\vf'(t) + f(\vf(t))\ \geq\ 0 \qquad{\rm a.e.}
\eqno{(\HH.12)}
$$
Note that at a point $t$ where $\vf$ is differentiable, if $\vf(t)=0$, then this
implies that $\vf'(t)\geq0$.  Thus (\HH.12) can be rewritten as 
$$
{-\vf'(t)  \over f(\vf(t))}\ \leq\ 1 \qquad{\rm a.e.}
\eqno{(\HH.13)}
$$
where the LHS equals $-\infty$ at points where $\vf(t)=0$.  Therefore, for any
measurable set $B$ we have
$$
-\int_B {\vf'(t)  \over f(\vf(t))}\ \leq \ |B|.
\eqno{(\HH.14)}
$$
On the set $B^-$ where $\vf$ is differentiable and $\vf'(t)<0$,  the inequality (\HH.14) has content.
Otherwise the integrand ${-\vf'(t)  \over f(\vf(t))} \leq 0$.  

Choose $s_1$ and $s_0$ so that 
 $r_1<s_1<s_0< r_0$ and $0< \vf(s_0) < \vf(s_1)$.
We will show that 
$$
 \int_{\vf(s_0)}^{\vf(s_1)}   {dy \over  f(y)}    \ \leq \  r_0-r_1\qquad{\rm for\ all\ such\ \ } s_0 > s_1,
\eqno{(\HH.15)}
$$
Because of (\HH.11) the point $s_0$ with $\vf(s_0)>0$ can be chosen arbitrarily close to $r_0$.
Then taking the limit as $s_0$ increases to $r_0$ proves that
$$
 \int_{0}^{\vf(s_1)}   {dy \over  f(y)}    \ \leq \  r_0-r_1 \ < \ \infty.
$$

It remains to prove (\HH.15).  Let $N(\vf\bigr|_A, y)$ denote the cardinality of
$\{ t\in A: \vf(t)=y\}$.  Set $A=[s_1,s_0]$, and let $V_A(\vf)$ denote the total variation of 
$\vf$ on $A$.  Since $\vf$ is absolutely continuous, we have, by Theorem 2.10.13 (p. 177)
in [\F], that 
$$
V_A(\vf)\ \  {\rm is\ finite,\  and\ \ } V_A(\vf)\ =\ \int N\left( \vf\bigr|_A, y   \right) \, dy
\eqno{(\HH.16)}
$$
Now set
$$
f_\e(y) \ \equiv\ \max\{f(y),\e\} \quad {\rm where}\ \ \e>0.
$$
Then 
$$
\int {1\over f_\e(y)} N\left( \vf\bigr|_A, y  \right) \,dy \ \leq \ { 1 \over \e } V_A(\vf)\ \ <\ \infty.
$$
Hence, the second half of Theorem 3.2.6 (p.245) in [\F] applies to yield
$$
 \int_{\vf(s_1)}^{\vf(s_0)} {1\over f_\e(y)} \, dy  
\ =\     -\int_{s_1}^{s_0}  { \vf'(t) \over f_\e(\vf(t))}\,dt.
\eqno{(\HH.17)}
$$

Since ${1\over f_\e(y)}\leq {1\over f(y)}$ on the set $B^-$ where $\vf$ is differentiable
and  $\vf'(t)<0$, we have
$$
\int_{B^-} {-\vf'(t) \, dt  \over  f_\e(\vf(t))} \ \leq \int_{B^-} {-\vf'(t) \, dt  \over  f(\vf(t))} 
\ \leq \ |B^-|\ \leq \ r_0 - r_1
\eqno{(\HH.18)}
$$
by (\HH.14).  Combining (\HH.17) and  (\HH.18) 
proves that 
$$
 \int_{\vf(s_1)}^{\vf(s_0)} {dy\over f_\e(y)} \ \leq \ r_0-r_1,
$$
since $\int_{\sim B^-}  {- \vf'(t)\, dt \over f_\e(\vf(t))} \leq 0$.
By the Monotone Convergence Theorem this proves (\HH.15).\qed

\Remark{\HH.7}  In the proof of Proposition \HH.5, the fact that
 $f$ is increasing was only used in Lemma \HH.6.  Therefore,
 if a subequation $E$ is defined by an  upper semi-continuous function
  $f:[0,\infty)\to [0,\infty]$   with $f(0)=0$ and $f(y)>0$ for $y>0$, then we 
  have that: 
  \medskip
  
  \centerline{\sl
  $ \int_{0^+} {dy\over f(y)} \ =\ \infty\ \Rightarrow$  the  (SMP)  holds}
  
  \centerline{\sl
for  all $E$-subharmonic functions $\psi$ for which $\psi'$ is absolutely continuous}.

 \vfill\eject
%\vskip.3in

%%%%%%%%%%%%%%%%%%%%%%%%%%%%%%%%%%%%%%%%%% 
%%%%%%%%%%%%%%%%%%%%%%%%%%%%%%%%%%%%%%%%%% 
%%%%%%%%%%%%%%%%%%%%%%%%%%%%%%%%%%%%%%%%%% 
%%%%%%%%%%%%%%%%%%%%%%%%%%%%%%%%%%%%%%%%%% 
%%%%%%%%%%%%%%%%%%%%%%%%%%%%%%%%%%%%%%%%%% 

\centerline{\headfont \II.\  Radial (Harmonic) Counterexamples to the (SMP)}.
 \smallskip

In this section we give the proof of Part (b) of Theorem \CC.6 by constructing
a radial counterexample to the (SMP) for $F$. 
Let $f$ denote the restriction of $\uf$ to $[0,\infty)$, where $\uf$ is the 
(smaller) characteristic function (Definition \CC.4) of the given borderline
subequation $F$.  Then
$$
f:[0,\infty) \to  [0, \infty]\ \ {\rm is\  upper\ semicontinuous, \   increasing \ and \  }f(0)=0.
\eqno{(\II.1)}
$$ 
More precisely we prove the following.

\Theorem{\II.1} {\sl
Suppose that $F$ is a borderline subequation with $f$ as described above.
If $\int_{0^+} {dy\over f(y)} <\infty$, then there exists a radially increasing $F$-subharmonic
function $u(x) = \psi(|x|)$ on $|x|>1$ where $\psi$ is of class $C^{1,1}$ on $(1,\infty)$
and satisfies}
$$
\psi(t) \ <\ m \ \ {\sl for}\ \ 1<t<t_0\quad{\sl and}\quad
\psi(t) \ =\ m  \ \ {\sl for}\ \  t \geq  t_0.
\eqno{(\II.2)}
$$ 

By Theorem \FF.3, it suffices to construct an increasing $C^{1,1}$-function which
is $R^{\ua}_f$-subharmonic and satisfies (\II.2).

In order to explicate the proof we will use the ``almost-everywhere theorem'' for
quasi-convex functions, which holds for the most general  possible subequations 
$F$.  This AE Theorem states that for a quasi-convex function $u$
$$
{\rm If\ } u \ {\rm has\  its\ 2-jet\  in\ } F \ {\rm a.e., \ then\ } u \ {\rm is\  } F\, {\rm subharmonic},
\eqno{(\II.3)}
$$ 
and was established in [\AETHM].  We will also make use of the fact (cf. [\HIR],  [\EB], or [\AETHM]) that
$$
 u \ {\rm is\  of\  class \   } C^{1,1}
 \qquad\iff\qquad
 u  \ {\rm and \ } -u \ {\rm are \ quasiconvex.  } 
 \eqno{(\II.4)}
$$

\noindent
{\bf Proof of Theorem \II.1.}  We start by solving the constant coefficient subequation
$E$ on $\bbr$ defined by
$$
E\ :\quad  a+f(p)\ \geq \ 0\and p\ \geq \ 0,
\eqno{(\II.5)}
$$ 
which is simpler than $R_f^\ua$.

\Lemma{\II.2}  {\sl
 If 
$
\int_{0^+} {dy\over f({y} )}<\infty
$,
then there exists an $E$-subharmonic function $\vf(s)$ of class
$C^{1,1}$ on $(0,\infty)$ with 
$$
\vf(s)<m \ \ {\rm strictly\ increasing \  on\ } (0,s_0)
\and
\vf(s) \equiv m \ \ {\rm on\ } [s_0, \infty)
$$
}

\pf
Set $s(y) = \int_0^y {dy\over f(y)}$ for $y\geq 0$.  
  For $0\leq y_1 < y_2\leq y_0$ we have
$$
 {y_2-y_1 \over  f(y_2)} \ \ \leq\ \ \int_{y_1}^{y_2} {d\,t\over f(t)}\  \ =\ \ s_2-s_1.
\eqno{(\II.6)}
$$ 
Therefore,  this function  $s(y)$ is strictly 
increasing until $f$ equals $+\infty$
(and is constant afterwards).  In particular, it is a homeomorphism from $[0,y_0]$ 
to $[0,s_0]$ for some $y_0>0$ with $s_0=s(y_0)<\infty$.  
Let $y(s)$ denote the inverse, which is also strictly increasing
with $y(0)=0$.  The inequality (\II.6) implies that  $y(s)$ is  Lipschitz 
on $[0, s_0]$ with Lipschitz constant $f(y_0)$, since $f(y_2)\leq f(y_0)$
if $y_2\leq y_0$.

Taking $y_1 =0$, $y_2=y(s)$ yields $y(s) \leq sf(y(s))$ which implies that
$y$ is differentiable  from the right at $s=0$ with $y'(0)=0$.  Moreover, since $y(s)$ is Lipschitz,
it is differentiable a.e. and 
$$
y'(s) = f(y(s))\quad a.e.
\eqno{(\II.7)}
$$

Fix $m$ and consider the function $\vf(s)$ defined on $(0, \infty)$ by $\vf(s_0)=m$ and 
$$
\vf'(s)\ \equiv\ 
\cases
{
y(s_0-s) \qquad  {\rm if}\ \  0  < s \leq  s_0  \cr
\qquad 0 \qquad\ \ \ \    {\rm if}\ \qquad \ s\geq s_0.
}
$$
Since $\vf'(s)$ is continuous and strictly decreasing to zero on $(0, s_0]$, $\vf(s)$ must be strictly
increasing to $m$ on  $(0,s_0]$ and identically equal to $m$ afterwards.

Since $\vf$ is twice differentiable at $s=s_0$, with $\vf'(s_0) = \vf''(s_0)=0$, the function
$\vf$ is class  $C^{1,1}$ on all of $(0,\infty)$.  Moreover, (\II.7) implies that
$$
\vf''(s) + f(\vf'(s))\ =\ 0\quad {\rm a.e.}  \ \   {\rm on}\ \ (0,\infty).
\eqno{(\II.8)}
$$ 
By (\II.4) and (\II.3) this implies that $\vf$  is $E$-subharmonic on $(0,\infty)$.
\qed
\medskip

We will use Lemma \II.2 applied to the subequation $E'$ defined by
$$
E'\ :\quad  a+ p + f(p)\ \geq \ 0\and p\ \geq \ 0,
\eqno{(\II.9)}
$$ 
rather than $E$. Now consider the  radial subequation $R_f^\ua$ on $(0,\infty)$
defined by
$$
R_f^\ua\ : \quad a+f\left( {p\over t} \right) \ \geq\ 0,
\and 
p\ \geq\ 0
\eqno{(\II.10)}
$$ 
which depends on the variable $t\in (0,\infty)$.

\Prop{\II.3}  {\sl
Suppose $\vf(s)$ is the $E'$-subharmonic function given by Lemma \II.2
applied to $E'$ rather than $E$.   Then the function $\psi (t)$
defined  on $(1,\infty)$ by
$$
\psi'(t)\ =\ t \vf'(\log\, t)
\and 
\psi(t_0)=m,
\eqno{(\II.11)}
$$ 
where $t_0= e^{s_0}$, is a $C^{1,1}$ subharmonic for $R_f^\ua$.  Moreover, 
$$
\eqalign
{
&\psi(t)\ \ {\rm is\ strictly\ increasing\ with\ }  \cr
\psi(t)\ <\ m\ \ {\rm on}\ \ &1<t<t_0 
\and
 \psi(t)\ \equiv\ m \ \ {\rm on}\ \  t_0  \leq t.
 }
\eqno{(\II.12)}
$$ 
}
\pf
That $\vf'$ is Lipschitz implies that $\psi'$ is Lipschitz.  Therefore $\psi$ is class $C^{1,1}$.
At a point of differentiability we have $\psi''(t) = \vf'(\log\, t) +  \vf''(\log\, t)$, and hence
$\psi''(t) +f({\psi'(t)\over t}) = \vf''(\log\, t) + \vf'(\log\,t) + f(\vf'(\log\, t))=0$.
Therefore $\psi(t)$ satisfies  (\II.10)  a.e.  
(Since $\vf'(s)$ is continuous and $>0$ on $(0,s_0)$, $\psi'(t)$ is also continuous and $>0$
on $(1,t_0)$. Thus $\psi$ is strictly increasing on $(1,t_0)$.)
Thus by (\II.4) and (\II.3), $\psi$ is $R^\ua_f$-subharmonic. 
The properties (\II.12) are straightforward.  \qed

\Remark {\II.4. ($F$-Harmonicity)}
The $F$-subharmonic function $u(x) = \psi(|x|)$ constructed in this section
is, in fact, $F$-harmonic if $F$ is invariant as in Definition \CC.7.  We leave it to
the reader to show that  $-\psi$ is ${\wt R}_f^\ua$-subharmonic and hence
$-u$ is $\ft$-subharmonic.  One can show that
$$
a+f(p) \ =\ 0, \ \ p\ \geq\ 0 
\qquad\Rightarrow\qquad
(p,a)\ \in\ \partial E,
\eqno{(\II.13)}
$$ 
but the converse is not  true if $f$ has a jump.
\medskip

\Ex{\II.5}
One of the simpler examples where Theorem \II.1 applies is the subequation $F$
defined by $\l_{\rm max}(A) \geq0$ and $\l_{\rm min}(A) + \sqrt{\l_{\rm max}(A)}  \geq0$ 
(See Example \JJ.1 below).  The characteristic function is $f(\l) = \sqrt\l$.
However, carrying out the construction of the $F$-harmonic counterexample provided in the
proof of Theorem \II.1 involves taking a complicated inverse.  To obtain more explicit harmonics
consider the subequation $F$ defined by $\l_{\rm max}(A) \geq0$ and
$\l_{\rm min}(A) + f({\l_{\rm max}(A)})  \geq0$ where
$$
f(\l) \ \equiv \ \sqrt \l \left( {4 R\over \sqrt{4R+\l} + \sqrt{\l}}    \right)
\ =\ \sqrt{\l^2+4R\l}  -\l.
$$
The characteristic function is $f(\l)$, and $\lim_{\l\to0} f(\l)/\sqrt\l = 2\sqrt R$ so that
$\int_{0^+}1/f <\infty$, and hence again Theorem \II.1 applies.

The increasing radial harmonics for this subequation $F$ on $\rn-\{0\}$ are very explicit:
$$
h(x)\ \equiv \ \cases{ -{R\over 3r}(r-|x|)^3 +k \qquad{\rm for}\ \ |x|\leq r \cr \ \cr
\qquad\qquad\ 0 \qquad\qquad\quad{\rm for}\ \ |x|\geq r }
\eqno{(\II.14)}
$$ 
A general version of this example is provided in   \JJ.13 in the next section.

 %\vfill\eject
\vskip.3in

%%%%%%%%%%%%%%%%%%%%%%%%%%%%%%%%%%%%%%%%%% 
%%%%%%%%%%%%%%%%%%%%%%%%%%%%%%%%%%%%%%%%%% 
%%%%%%%%%%%%%%%%%%%%%%%%%%%%%%%%%%%%%%%%%% 
%%%%%%%%%%%%%%%%%%%%%%%%%%%%%%%%%%%%%%%%%% 
%%%%%%%%%%%%%%%%%%%%%%%%%%%%%%%%%%%%%%%%%% 

\centerline{\headfont \JJ.\  Subequations with the Same Increasing Radial Subharmonics.}.
 \smallskip

In order to begin to understand examples and applications of the main Theorem A,
it is helpful to describe all the borderline invariant subequations with a given
 characteristic function $f$.
 
 \medskip
 \noindent
 {\bf Remark.}  By Theorem \GG.3, this problem is equivalent to describing all borderline invariant
 subequation with the same set of increasing radial subharmonics (or, equivalently, the same
set of increasing radial harmonics) satisfying
$$
  R_f^\ua \ \  :\ \ \ \ \psi'(t) \ \geq\ 0 \and \psi''(t)  + f\left({\psi'(t) \over t}\right) \ \geq\ 0
\eqno{(\JJ.1 )}
$$
on $(\a,\b) \ss(0, \infty)$. 
\medskip

We assume that an increasing upper semi-continuous function
$$
f: [0, \infty) \ \to\ [0, \infty)\qquad{\rm with}\ \ f(0)\ =\ 0
\eqno{(\JJ.2)}
$$
is given.  The problem is to determine all subequations (if any) with this
characteristic function $f$.

We start with the two examples that play a central role.  
Given $A\in\Symn$,  let $\l_1(A)\leq\cdots \leq \l_n(A)$ 
denote the ordered eigenvalues of $A$.
In particular, the minimum and maximum eigenvalues are $\l_{\rm min}(A) =\l_1(A)$
and $\l_{\rm max}(A) =\l_n(A)$ respectively.
Recall the monotonicity $\l_k(A+P) \geq \l_k(A)$ for $P\in \cp$.

\Ex{\JJ.1. (The $f$-Min/Max Subequation)}
$$
F_f^{\rm min/max} \ \equiv\ \{ A : \l_{\rm max}(A) \geq 0 \ \ {\rm and}\ \ 
\l_{\rm min}(A) +f(\l_{\rm max}(A) ) \geq0\}
$$

\Ex{\JJ.2. (The $f$-Min/2 Subequation)}
$$
F_f^{\rm min/2} \ \equiv\ \{ A : \l_{2}(A) \geq 0 \ \ {\rm and}\ \ 
\l_{\rm min}(A) +f(\l_{2}(A) ) \geq0\}
$$

\Prop{\JJ.3} {\sl
The sets $F_f^{\rm min/max}$ and $F_f^{\rm min/2}$ are subequations which
are borderline and O$_n$-invariant.  Moreover, for both subequations, the characteristic function
restricted to $[0,\infty)$  equals  $f$.
}

\pf 
Since $f$  is upper semi-continuous, both sets are closed.
Since $f$ is  increasing, positivity follows from the $\cp$-monotonicity of the ordered
eigenvalues. To prove these subequations are  borderline,
suppose  $A$ lies in the larger subequation  $F_f^{\rm min/max}$ 
and  $A\in -\cp$, i.e.,  $\l_{\rm max}(A) \leq 0$. Then  $\l_{\rm max}(A) = 0$
and since $f(0)=0$,  $\l_{\rm min}(A) = 0$.  Hence, $A=0$.
Invariance follows because  the ordered eigenvalues
themselves are O$_n$-invariant.   

To complete the proof we compute the full radial profiles (not just the
increasing part).

\centerline
{The subequation $F_f^{\rm min/max}$  has radial profile}
$$
\{(\l,\mu) : \l\geq 0\ \ {\rm and}\ \ \mu+ f(\l)\geq0\} \cup
 \{(\l,\mu) : \mu\geq 0\ \ {\rm and}\ \ \l+ f(\mu)\geq0\}.
\eqno{(\JJ.3)}
$$ 

\centerline{For $n\geq3$, the subequation $F_f^{\rm min/2}$ has radial profile}
$$
 \{(\l,\mu) : \l\geq 0\ \ {\rm and}\ \ \mu+ f(\l)\geq0\}.
\eqno{(\JJ.4)}
$$  

We see this as follows.
Note that the radial profile of $F_f^{\rm min/max}$ is symmetric about the diagonal.
Recall that if $A\equiv \l P_{e^\perp} + \mu P_e$ belongs to any borderline subequation,
then either $\l\geq0$ or $\mu\geq0$.

For (\JJ.3), suppose $A\equiv \l P_{e^\perp} + \mu P_{e} \in  F_f^{\rm min/max}$.
If $\l\geq \mu$, then $\l=\l_{\rm max}\geq0$ and $\mu=\l_{\rm min}$ satisfies
$\l_{\rm min} + f(\l_{\rm max})\geq0$.
If $\mu\geq \l$, then $\mu=\l_{\rm max}\geq0$ and $\l=\l_{\rm min}$ satisfies
$\l_{\rm min} + f(\l_{\rm max})\geq0$.

For (\JJ.4), suppose $A\equiv \l P_{e^\perp} + \mu P_{e} \in  F_f^{\rm min/2}$.
Since $n\geq3$, $\l =\l_2\geq0$ and either $\mu=\l_1$ or $\mu>\l$.
In either case $\l\geq0$ and
$\mu+f(\l)\geq0$.\qed

\Cor{\JJ.4} {\sl
Both $F_f^{\rm min/max}$ and $F_f^{\rm min/2}$ have their increasing radial subharmonics
$u(x) =\psi(|x|)$ determined by the subequation $R^{\uparrow}_f$ defined in (\JJ.1).}

\medskip

The subequations $F_f^{\rm min/max}$ and $F_f^{\rm min/2}$ are of central importance because they 
are the largest and smallest possible under our  invariance  hypothesis (\CC.9)  on $F$:
$$
\l P_{e^\perp} +\mu P_e \in F\ \ {\rm for\ some\ } e\neq0 
\qquad\Rightarrow\qquad
\l P_{e^\perp} +\mu P_e \in F\ \ {\rm for\ all\ } e\neq0.
\eqno{(\CC.9)}
$$

\Theorem {\JJ.5} { \sl 
Suppose $F$ is invariant.
Then $F$ has characteristic function $f$, or equivalently, 
the radial increasing subharmonics 
$u(x) =\psi(|x|)$ for $F$ are determined by $R^{\uparrow}_f$  as in (\JJ.1),
if and only if }
$$
F_f^{\rm min/2} \ \ss\ F\ \ss\ F_f^{\rm min/max}.
\eqno{(\JJ.5)}
$$  
\pf
Each $A\in\Symn$ can be written as  a sum $A= \l_1 P_{e_1} + \cdots + \l_n  P_{e_n}$
using the ordered eigenvalues of $A$. 
Set $B_0 \equiv \l_1 P_{e_1} + \l_2 P_{e_1^{\perp}}$ and 
$B_1 \equiv \l_1 P_{e_1} + \l_n P_{e_1^{\perp}}$, and note that  $B_0\leq A\leq B_1$.

If $A\in F_f^{\rm min/2}$, then $\l_2\geq0$ and 
$\l_1 + f(\l_2)\geq0$.  Thus $B_0 \equiv \l_1 P_{e_1} + \l_2 P_{e_1^{\perp}} \in F_f^{\rm min/2}$.  
Since $F_f^{\rm min/2}$ and $F$ have the same radial profile
in the half-plane $\{\l\geq0\}$ by (\JJ.4), we conclude that $B_0\in F$.
However, $B_0\leq A$ proving that $A\in F$.

For the other inclusion, pick $A\in F$.  Since $F \ss\cpt$ we have $\l_{\rm max}\geq0$.
Now $A\leq B_1$  implies $B_1\in F$.
By the invariance  hypothesis and  (\JJ.3), $F$ and $F_f^{\rm min/max}$ have the same
same radial profile in the half-plane $\{\l\geq0\}$. Therefore, $B_1 \in F_f^{\rm min/max}$,
i.e., $\l_n\geq0$ and $\l_1+f(\l_n)\geq0$. 
This implies by definition that $A\in F_f^{\rm min/max}$. \qed

\Remark{\JJ.6} Theorem \JJ.5 can be used to construct vast numbers of invariant 
borderline subequations which satisfy the (SMP), or, if one prefers, which do not satisfy the (SMP).

\Remark{\JJ.7} Dropping the invariance assumption (\CC.9), the proof of Theorem \JJ.5
shows that for any borderline subequation $F$ with characteristic functions $\uf$ and $\of$
$$
F_{\uf}^{\rm min/2}\ \ss\ F\ \ss\ F_{\of}^{\rm min/max}.
\eqno{(\JJ.5)'}
$$

The subequation $\Fh$ defined by (\HH.1) as the O$_n$-orbit of $F$, has characteristic
function $\of$.  It satisfies 
$$
F \ \ss\ \Fh\ \ss\ F_{\of}^{\rm min/max}.
\eqno{(\JJ.6)}
$$
and is the smallest  O$_n$-invariant subequation containing $F$.

\bigskip

%%%%%%%%%%%%%%%%%%%%%%%%%%%%%%%%%%%%%%%%%% 
%%%%%%%%%%%%%%%%%%%%%%%%%%%%%%%%%%%%%%%%%% 
%%%%%%%%%%%%%%%%%%%%%%%%%%%%%%%%%%%%%%%%%% 
%%%%%%%%%%%%%%%%%%%%%%%%%%%%%%%%%%%%%%%%%% 
%%%%%%%%%%%%%%%%%%%%%%%%%%%%%%%%%%%%%%%%%% 

\centerline{\bf  Explicit  Borderline Examples.}.
 \smallskip

It is natural to look for the largest possible subequations which  satisfy the (SMP).  
Because of Theorem \CC.6(a) and Theorem \JJ.5
these are max-min subequations $F$ whose characteristic function $f$  is as large as 
possible subject to the condition $\int_{0^+}{dy\over f(y)} = \infty$.  Since this only depends
on the behavior of the germ at $0^+$ of $f$, we can also localize $F$ at the origin in $\Symn$
by replacing $f$ on $[\e,\infty)$ by the function $\equiv +\infty$, which yields a larger subequation.
We present three examples where the (SMP) holds, and two where the  (SMP) fails.

If $F$ is a cone, then by (\CC.15) $\of(\l) = \overline\a \l$, and the corresponding min/max-subequation
$F^{\rm min/max}_{\of}$ containing $F$ is given as follows.

\Ex {\JJ.8. (Min/Max Cones)} \ ($0<\a<\infty$)
$$
   f(y) = \a y \qquad   \left ({\rm and\ hence\ } \int_{0^+} {1\over f} \ =\ \infty\right)
\leqno{(a)}
$$
The borderline O$_n$-invariant cone subequation
$$
   \cp_{\a}^{\rm min/max} \ :\quad   \l_{\rm min}(A) + \a \l_{\rm max}(A)\  \geq \ 0 \qquad{\rm satisfies\ the\ (SMP)}.
\leqno{(b)}
$$
Thus all subequations $F$ which are contained in $\cp_{\a}^{\rm min/max}$ for some $\a>0$ satisfy the (SMP).
 The increasing radial harmonics $\psi(|x|)$ are important classical functions given (with $p\equiv \a+1$) by:
$$
\eqalign
{
&  \psi(t) \ =\ a t^{2-p} +b\quad   {\rm with}\ \ a>0 \ \ {\rm if}\ \  1\leq p <2 \cr
&     \psi(t) \ =\ a \log t +b\quad   {\rm with}\ \ a>0 \ \ {\rm if}\ \  p=2 \cr 
&  \psi(t) \ =\ - {a \over  t^{p-2} }+b\quad   {\rm with}\ \ a>0 \ \ {\rm if}\ \  p>2 \cr
}
\eqno{(\JJ.7)}
$$

We leave the computation of $M_F$ for this $F$ as an open problem.

This example can be localized.

\Ex {\JJ.9. (Localizing  the Min/Max Cone)} \ ($0<\a<\infty$ and $\e>0$)
$$
   f(y) \ \equiv\ 
   \cases
   {
   \a y \qquad {\rm if} \ \ 0\leq y<\e  \cr
 +  \infty \qquad {\rm if} \ \ \e\leq y
   }
\leqno{(a)}
$$\medskip\noindent
(b)\ \ The borderline O$_n$-invariant cone subequation
$$
   \cp_{\a, {\rm loc}}^{\rm min/max}\ :\quad   {\rm Either \ \ } \l_{\rm min}(A) \ \geq \ \e \ \ \ {\rm or}\ \ \ 
    \l_{\rm min}(A) + \a \l_{\rm max}(A)\  \geq \ 0 % \qquad{\rm satisfies\ the\ (SMP)}.
% \leqno{(b)}
$$
satisfies the (SMP).  Thus
$$
{\rm All \ subequations}\ F\ {\rm which\  are\ contained\ in\ }   \cp_{\a, {\rm loc}}^{\rm min/max}
\ {\rm for\ some\ }\a \ {\rm satisfy\ the\ (SMP).}
\eqno{(\JJ.8)}
$$

\Remark{\JJ.10. (The Barles and Busca Hopf Lemma 3.2 [\BaB])}
Under their hypothesis ``(F3b)'' they prove that the (SMP) holds.  This landmark
paper on comparison covers a wide range of subequations.  For the constant coefficient,
pure second-order subequations considered here their hypothesis can be restated as follows:
$$
\eqalign
{
&\qquad\qquad\forall\, \l> 0, \ \ \exists\, \mu,   \d\ >\ 0\ \ \ {\rm such\ that\ }\cr
&\ \ \ \ E\ \equiv\ \left\{t  \left (\l P_{e^\perp} - \mu P_e\right) : 0<t<\d\ \  {\rm and}\ \ |e|=1\right\}  \cr
&{\rm is\ contained\ in\ the \ complement \ of\ the\ subequation\ } F.
}
\eqno{(F3b)}
$$
Now the assertion that $E\ss ( \sim F)$ is equivalent to saying that $\of(y) < {\mu \over \l} y$
for all $0<y<\l\d$.
By  Theorem \JJ.5 and Remark \JJ.7 this proves that the
condition (F3b) is equivalent to 
$$
F\ \ss  \ \cp_{\a, {\rm loc}}^{\rm min/max}  \qquad {\rm for\ some\ \ } \a
$$
(take $\a={\mu  \over  \l}$).
Thus, the Hopf Lemma (3.2) in [\BaB], 
when restricted to subequations of the type considered here,
is equivalent to the corollary (\JJ.8) of Theorem \CC.6(a).

\Ex {\JJ.11. (A Localized Hopf Subequation)} \ ($0< k < \a<\infty$ and $0<\e\leq 1$)
$$
   f(y) =
   \cases
   {
     y\left(  \a  + k\log{1\over y} \right) \  \qquad {\rm if} \ \ 0\leq y<\e  \cr
 +  \infty \qquad \qquad \qquad \quad {\rm if} \ \ \e\leq y
   }
\leqno{(a)}
$$
is an upper semi-continuous increasing function with associated min/max subequation
$$
   H(\a) \ :\quad   {\rm Either \ \ } \l_{\rm min}(A) \ \geq \ \e \ \ \ {\rm or}\ \ \ 
    \l_{\rm min}(A) + \l_{\rm max}(A) \left( \a - k \, \log \,{  \l_{\rm max}(A) }\right)
    \  \geq \ 0.
 \leqno{(b)}
$$
This subequation satisfies the (SMP) since
$$
\int{dy\over y(\a-k\log y)} \ =\ -{1\over k} \log (\a-k\log y)
$$
which implies $\int_{0^+}{du\over f(y)} = \infty$.    The increasing radial harmonics satisfy
$$
\psi'(t) \ =\ \b t e^{ct^k} \quad {\rm where \ \ } \log \, \b \ =\ {1+\a\over k}
\eqno{(\JJ.9)}
$$
and $c$ is the constant of integration.  For example, if we  take $k=2$ and set  $c=-\b/2$, we see that 
$\psi'(t) = \b t e^{-\b t^2\over 2}$
integrates to 
$$
\psi(t) \ =\  e^{-\b R^2\over 2} - e^{-\b t^2\over 2}
\eqno{(\JJ.10)}
$$
which is the standard Hopf function (cf. [\CLN]).
Here $f(y) = y(2\log ({\b\over y})-1)$ for $y$ small.
\medskip

Obviously,
$$
 \cp_{\a}^{\rm max/min}  \ \ss\ \cp_{\a}^{\rm loc}\ \ss\ H_\a,
\eqno{(\JJ.11)}
$$
and for larger $\a$ each subequation is larger.
Our notation suppresses  the dependence of $\cp_{\a}^{\rm loc}$ on $\e$ and of 
$H_\a$ on $\e$ and $k$.

\Ex{\JJ.12. (The (SMP) Fails)} 
Let  $f: [0,\infty) \to  [0,\infty)$ be defined by
 $$
f(y) \ =\ N y^{ N-1\over N} \qquad (N\ >\ 1).
\eqno{(\JJ.12)}
$$
Since $\int {dy \over f(y)}  = y^{1\over N}$, we have $\int_{0^+} {dy \over f(y)} <\infty$.
Therefore the (SMP) fails for the corresponding min/max subequation
 $$
F\ :  \l_{\rm min}(A) + N \l_{\rm max}(A)^{N-1\over N}  \ \geq \ 0 \quad{\rm and }\quad \l_{\rm max}(A) >0.
\eqno{(\JJ.13)}
$$
The associated constant coefficient subequation
$$
E\  \ :\quad a\ +\ Np^{N-1\over N} \and p\ \geq\ 0
$$
has radial harmonics
$$
\psi(t) \ \equiv\ -{1\over N+1} (R-t)^{N+1} + k \quad{\rm for}\ \ t\leq R
\and
\psi(t) \ \equiv\ 0 \quad{\rm for}\ \ t\geq R,
\eqno{(\JJ.14)}
$$
but the radial harmonics for $R^\uparrow_f$ are more complicated.
Note that a better example (i.e.,  $f$ is smaller)
where (SMP) fails is $f(y) \equiv y (\log y)^2  <  y^\b$
($0<\b<1$ and $y$ small), since $\int{dy\over f(y)} = 1+ {1\over \log y} <1$ for $y>0$ small.

\Ex{\JJ.13. (The (SMP) Fails with Explicit Harmonics)}
However, fixing $R>0$ there exists a modification of (\JJ.13) of the form
$$
F \ \ : \quad  \l_{\rm min}(A) + N(\l_{\rm max}(A) g(\l_{\rm max}(A) ))^{N-1\over N} \ \geq\ 0
\and
\l_{\rm max}(A) \ \geq\ 0
\eqno{(\JJ.15)}
$$
with $g$ defined below so that
$F$ has the simple/explicit radial harmonics:
$$
\psi(|x|) \ \equiv \ -{r^2\over (N+1)R^2}(r-|x|)^{N+1} + k \ \ {\rm for } \ |x|\leq r
\and
\psi(|x|)\equiv 0 \ \ {\rm for }\ |x|\geq r
\eqno{(\JJ.16)}
$$
The proof is omitted.

The characteristic function for $F$ is 
$$
f(y)\ =\ N(y g(y))^{N-1\over N}.
\eqno{(\JJ.17)}
$$
The function $g(y)$ is defined to be the inverse of $y(t) \equiv (R-t)^N/t = \psi'(t)/t$
for $0\leq t\leq R$.
Since $y(t)$ is strictly decreasing from $\infty$ to $0$ on $[0,R]$, the function $g(y)$
is strictly decreasing from $R$ to $0$ on $[0, \infty]$.
One can show that $x={1\over N} f(y)$ has inverse 
$$y(x) = {x^{N\over N-1}  \over R-x^{N\over N-1} },
$$
and hence is strictly increasing on $[0,\infty]$ from $0$ to $NR^{N-1}$ insuring that $F$ is a subequation.

 \vfill\eject

% \vfill\eject
\vskip.3in

%%%%%%%%%%%%%%%%%%%%%%%%%%%%%%%%%%%%%%%%%% 
%%%%%%%%%%%%%%%%%%%%%%%%%%%%%%%%%%%%%%%%%% 
%%%%%%%%%%%%%%%%%%%%%%%%%%%%%%%%%%%%%%%%%% 
%%%%%%%%%%%%%%%%%%%%%%%%%%%%%%%%%%%%%%%%%% 
%%%%%%%%%%%%%%%%%%%%%%%%%%%%%%%%%%%%%%%%%% 

\centerline{\headfont \KK.\   Strong Comparison and Monotonicity.}.
 \medskip

By the {\bf strong comparison principle} for a {\bf subequation} $F$ we mean the following.
$$
{\rm If}\ u\in F(\ob)\ {\rm and}\ v\in \ft(\ob), \ {\rm then\ the \ (SMP) \ holds\ for\ } u+v \ {\rm on} \ \ob.
\eqno{(SC)}
$$
This is, of course, immediate if $u+v$ is $G$-subharmonic for some subequation $G$
for which the strong maximum principle holds. 
In this section we address the question of when such a $G$ exists.
The geometric point of view is, we think, an advantage here.  This question is reduced to algebra
by the following.

\Theorem {\KK.1. (Addition)}
{\sl If  three subequations satisfy
$$
F+H\ \ss\ G,
$$
then}
$$
F(\ob)  +  H(\ob)  \ \ \ss\ \  G(\ob).
$$

\medskip
\noindent
{\bf Remark.}
This result is immediate from  sup-convolution and either of the 
classical Jensen or Slodkowsky Lemmas (which are in  a
strong sense equivalent, cf. [\AETHM]).
It is referred to as ``Transitivity of inequalities in the viscosity sense'' on
page 745 of  [\AS], and is proved in the  book [\CC] of Caffarelli-Cabr\'e (Proposition 2.9 )
 in the case where $F$ and $H$ are uniformly elliptic.
See also classic works of Crandall and Crandall, Ishii and Lions [\C], [\CIL].
\medskip

Thus the (SC) question is reduced to asking when is 
$F+\ft$ contained in $G$, where $G$ satisfies the (SMP).

\medskip

Using the fact (\BB.6) that $\wt{F+A} = \ft -A$ one can show that for any two subequations $F$ and $G$
$$
F+\ft \ \ss\  G \qquad\iff\qquad F + \wt G \ \ss\  F.
%\eqno{(\KK.1)}
$$
Rewriting this with $\wt G$ replaced by $M$ gives
$$
F + M \ \ss\  F  \qquad\iff\qquad  F+\ft \ \ss\  \wt M.
\eqno{(\KK.1)}
$$

A subequation $M$ satisfying $F+M\ss F$ will be called a {\bf monotonicity subequation for $F$}.
It is easy to show that $M$ is a monotonicity subequation for $F$ if and only if  $M$ is
monotonicity subequation for $\ft$.  (See (5) below.)

\Theorem{\KK.2. (Strong Comparison)}
{\sl
Suppose  that $M$ is a monotonicity subequation for $F$.   Then}
$$
{\rm (SMP) \ \ for \ \ } \wt M
 \qquad\Rightarrow\qquad
{\rm (SC) \ \ for \ \ } F 
\eqno{(\KK.2)}
$$
\pf
By   (\KK.1) and Theorem \KK.1, $F+M\ss F \ \Rightarrow\  F+\ft\ss \wt M \ \Rightarrow\  F(\ob)+\ft(\ob)\ss \wt M (\ob)$.\qed

\vfill\eject
\centerline{\bf The Largest Monotonicity Subequation for $F$}
\medskip

Increasing the size of a subequation $M$ satisfying $F+M\ss F$ decreases the size of 
$G=\wt M$, thereby increasing the liklyhood that  $G=\wt M$ satisfies the (SMP).  Hence,
it is natural  to look for the largest subequation $M$ satisfying $F+M\ss F$. 
It is somewhat surprising that there is such a subequation.
We define the  {\bf monotonicity subequation for}  $F$ to be the set
$$
M_F \ \equiv\  \{A\in\Symn : F+A\ss F\}.
\eqno{(\KK.3)}
$$
We leave the following facts as an exercise.
$$
\eqalign
{
&(1)\ \ M_F \ {\rm is\ a\ subequation}, \qquad (2) \ \ 0\in \partial M_F \ \ {\rm and}\ \ \wt{M}_F \ \ss\ \cpt,   \cr
&(3) \ \ M_F \ {\rm is\  its\  own\ monotonicity\  subequation,\  in \  particular\ } M_F \ {\rm is \ additive,}  \cr
&(4) \ \ {\rm If }\ M_F \ {\rm is\  a\  cone,\  then\ } M_F \ {\rm is\ a \ convex\ cone\   subequation,}  \cr
&(5) \ \ M_{\ft} \ =\ M_F, \  \ {\rm in\ fact, \ for\ any}\ A, \ \  \  F+A\ \ss\ F \quad\iff\quad \ft+A \ \ss\ \ft,  \cr
&(6) \ \  \Int \wt{M}_F\ \ss\ F+\Int \ft \ \ss\ F+\ft\ \ss\  \wt{M}_F,\ \ \ {\rm and\  hence}\ \ 
\wt M_F \ =\ {\rm Cl}(F+\ft) \cr
}
\eqno{(\KK.4)}
$$

\Def{\KK.3}  A subequation $M$ such that $0\in M$ and $M$ is additive, i.e., $M+M\ss M$,
will be called a {\bf monotonicity subequation}.

$$
\eqalign
{
&\qquad M\   {\rm is\ a \  monotonicity \ subequation} \quad\iff \cr
& M\ =\ M_F \ \  {\rm for\ some\ subequation\ } F.
}
\eqno{(\KK.5)}
$$
\pf 
In fact  $M$ is its own monotonicity subequation, because if $M+A\ss M$, then 
$0\in M \Rightarrow A\in M$. \qed
\medskip

Note that $M_F$ is maximal, that is, it contains every monotonicity subequation for $F$.
Consequently, Theorem \KK.2 could be restated equivalently as follows.

\Theorem{\KK.2$'$}
$$
\wt M_F \ \ {\rm satisfies\  the \ (SMP) }
\quad\Rightarrow \quad 
F\ \ {\rm satisfies  \ (SC)}
\eqno{(\KK.2)'}
$$

For most subequations $F$, even when $F$ is not a cone, $M_F$ is a cone.
These subequations $F$ will be referred to as {\bf normal subequations}.
If $F$ is normal, then in fact, by (4) above, $M_F$ is a convex cone.  
The verious criteria in Theorem \CC.9 apply to $\wt M_F$.  In addition,
uniform ellipticity can be added to the list since $M_F$ is a convex cone.

\Prop{\KK.4} 
{\sl   Suppose $F$ is a normal subequation (i.e., $M_F$ is a cone).  Then}
$$
\eqalign
{
{\sl (SMP)\ holds \  for \ } \ \wt M_F   \qquad  &\iff\quad -P_e  \notin \wt M_F \ \forall\,e\neq 0
  \qquad  \iff\qquad P_e \in\Int M_F \ \forall\,e\neq 0     \cr
 &\iff\qquad M_F \ \ {\sl is\ a\  convex\  conical \ neighborhood \ of\ \ } \cp \cr
  &\iff\qquad F \ \ {\sl is\ uniformly\ elliptic}.
}
$$
\pf
The first equivalence is just part (b) of Theorem \CC.9.
The second follows from the definition of the dual of $M_F$.  The third
follows since $M_F$ is a convex cone.  The last follows from Lemma B.1
in Appendix B.\qed
\medskip

We note that since $M_F$ is maximal, there is the possibility that the reverse implication
in (\KK.2)$'$ holds.   We leave this as an open question even  in the case where $F$ is a cone.
However, the following is a partial answer in this case.

\Prop {\KK.5} {\sl
Suppose that $F$ is a normal subequation and $F+\ft = \wt M_F$. Then
\medskip
\centerline {
(SC) holds for $F \quad\iff\quad$
the (SMP) holds for $\wt M_F$.
}
}

\pf
Suppose that the (SMP) fails for $\wt M_F$.  Then by
 Proposition \KK.4 we have $-P_e \in \wt M_F$ for some $e$. By the hypothesis that $F+\ft = \wt M_F$, we have
$$
-P_e\ =\ Q + \wt Q \qquad{\rm with}\ \ Q\in F \ \ {\rm and}\ \ \wt Q\in \ft.
\eqno{(\KK.6)}
$$
Let $w(x) \equiv \half {\bra {e}x}^2, \ u(x)  \equiv \half \bra {Qx}x$, and $v(x)  \equiv \half \bra {\wt Qx}x$
denote the corresponding quadratic functions.  Then $w=u+v$, $u\in F(\rn), v\in \ft(\rn)$
but the (SMP) fails for $w$.  Hence, $u$ and $v$ provide a counterexample to (SC) for $F$.
\qed
\medskip

The following corollary probably comes as no surprise.

\Cor{\KK.6}
{\sl  Suppose that $F$ is a convex cone subequation.  Then
\medskip
\centerline
{
(SC) holds for $F$ \qquad$\iff$\qquad $F$ is uniformly elliptic.
}
}
\pf
If $A, B \in F$, then $\half(A+B) \in F$ by convexity, and since $F$ is a cone, this proves
$F+F\ss F$  which implies $F\ss M_F$.  Since $0\in F$, $M_F = 0+M_F\ss F+M_F\ss F$.
Thus, $M_F=F$, and so $\wt M_F = \ft$.  
By (5) above, $\ft+ F \ss \ft  = \wt M_F$, while $\wt M_F = \ft \ss \ft+0 \ss \ft+F$.
This proves that $\ft+F = \wt M_F = \ft$ so that Proposition \KK.5 applies.
Finally, by Proposition \KK.4 the (SMP) holds for $\ft =\wt M_F \iff F$ is uniformly elliptic.\qed

\medskip

The hypothesis $F+\ft = \wt M_F$  in Proposition \KK.5
can be analyzed further because  the set $F+\ft$ can be explicitly
computed.  For this we introduce the  {\bf strict monotonicity set $S_F$ for $F$}
$$
S_F\ \equiv\ \{A\in\Symn : F+A \ss \Int F\},
\eqno{(\KK.7)}
$$
along with a secondary notion for the dual, this time for an arbitrary subset $G$, namely,
$$
G^* \ \equiv \ -(\sim G) \ =\ \sim(-G).
\eqno{(\KK.8)}
$$

\Remark{\KK.7}  Although we restrict  attention in this paper to  subequations $F\ss \Symn$,
it is worth noting that the next result holds for an arbitrary subequation $F\ss J^2(X)$ on a manifold $X$.

\Lemma{\KK.8} {\sl  For any subequation $F$}
$$
F + \ft\ =\ S_F^*.
$$
\pf
Note that 
$E\in S_F^* \iff -E\notin S_F  \iff   \exists A\in F$ such that $-B=A-E \notin \Int F
\iff E=A+B$ for some $A\in F$ and $B\in\ft$.\qed
\medskip

\Remark {\KK.9. (Reformulating $F+\ft=\wt M_F$)}  {\sl
For any subequation $F$ the following are equivalent statements.
\medskip

\centerline
{
(1) \ \ \ $\wt M_F \ \ss\ F+\ft$\qquad
(2) \ \ \ $\wt M_F \ \ss\ S_F^*$ \qquad
(3)\ \ \  $S_F \ \ss\ \Int M_F$
}
\medskip
\centerline
{
(4) $S_F \cap \partial M_F\ =\ \emptyset.$
}
\medskip\noindent
The reverse containments in (1), (2) and (3) are always true.  Thus if any of 
(1) through (4) are true, then equality holds in (1) through (3).
}
\medskip

Assertions (1) and (2) are equivalent since $S_F^* = F+\ft$
Assertions (2) and (3) are equivalent since $(\wt M_F)^* =\Int M_F$, 
$(S_F^*)^* = S_F$, and taking the secondary dual $(\cdot)^*$ reverses containments.
Assertions (1) and (2) are equivalent since $S_F \ss M_F$. \qed

\Ex{\KK.10. ($F+\ft = \wt M_F$ is true)}
Let $F$ be a convex cone subequation.
We saw, in the proof of Corollary \KK.6, that  $F+\ft = \wt M_F$.
Here we give a second proof involving $S_F$.
Note that $S_F = \Int M_F= \Int F$,
since any $B\in \partial M_F = \partial F$ cannot be in $S_F$ (because $B\in S_F$ 
would imply $0+B\in\Int F$ by definition).  Therefore $S_F^*=\ft =\wt M_F$.
Now apply Lemma \KK.8.

\medskip\noindent
{\bf Example \KK.11. ($F+\ft=\wt M_F$ is false).}
The simplest example is the Monge-Amp\`ere equation $F: \det A\geq1, A> 0$.
Here $M_F = \cp$, but $S_F = \cp-\{0\}$ is larger than $\Int M_F$. Thus, 
$S_F^*$ is smaller than $\wt M_F$.  In fact,
$S_F^* = (\Int \cpt)\cup \{0\} = F+\ft$ does not contain $-P_e$
 for any $e\neq 0$, but $-P_e\in \wt M_F$.
\medskip

This example does not contradict the equivalence of (SC) for $F$ and the (SMP) for $\wt M_F$
since both conditions fail in this case.  To see that (SC) fails for the Monge-Amp\`ere equation
$F$, one employs the classical Pogorelov harmonics
$$
h_\a(t,x) \ \equiv\ {|x|^{2-{2\over n}}  \over f_\a(t)^{1-{2\over n}}}
\qquad{\rm where} \ \ \ f_\a'' +f_\a^{n-1} \ =\ 0, \ \ {\rm and}\ \ f_\a(0)=\a>0.
$$
Then near $t=0$, $h_{2\a}-h_{\a} \leq0$ attains the maximum value zero on the 
$t$-axis.

\medskip

This section  leads to the question:  are there examples where 
$M_F$ is not a cone (i.e., $F$ is not normal)? The answer 
is yes.  This involves some  intriguing new subequations discussed in the next section.

 %\vfill\eject
\vskip.3in

%%%%%%%%%%%%%%%%%%%%%%%%%%%%%%%%%%%%%%%%%% 
%%%%%%%%%%%%%%%%%%%%%%%%%%%%%%%%%%%%%%%%%% 
%%%%%%%%%%%%%%%%%%%%%%%%%%%%%%%%%%%%%%%%%% 
%%%%%%%%%%%%%%%%%%%%%%%%%%%%%%%%%%%%%%%%%% 
%%%%%%%%%%%%%%%%%%%%%%%%%%%%%%%%%%%%%%%%%% 

\centerline{\headfont \LL.\   Examples of  Exotic Monotonicity Subequations which are not Cones}.
 \medskip

The examples will be constructed as follows.

\Def{\LL.1}  Suppose $g:[0,\infty)\to \bbr$ is a continuous decreasing function with 
$g(0)=0$ and $g(x)<0$ for $x>0$.  Set
$$
M^g   \ \equiv\ \{A : \  \tr A\geq0  \ \ {\rm and }\ \   \l_{\rm min}(A) \geq g( \tr A) \}
\eqno{(\LL.1)}
$$

\Prop{\LL.2}
{\sl  
$M^g$ is a subequation which is orthogonally invariant with
$$
M^g \cap \{\tr A=0\} \ =\ \{0\}
\and
\cp - \{0\} \ \ss\ \Int M^g.
\eqno{(\LL.2)}
$$
}
\pf  Since $g$ is continuous, $M^g$ is a closed set.  Recall that
$$
 \l_{\rm min}(A+B) \geq   \l_{\rm min}(A) + \l_{\rm min}(B).
\eqno{(\LL.3)}
$$
This combined with the fact that $g$ is decreasing easily implies that
positivity (P) holds for $M^g$.  Obviously $M^g$ is O$_n$-invariant.

If $A\in M^g$ and $\tr A=0$, then since $g(0)=0$, the minimum eigenvalue
$ \l_{\rm min}(A) \geq g(0)=0$.  But then $\tr A=0$ implies $A=0$.

If $P\geq0$ and $P\neq0$, then $\tr P>0$. Since $x>0$ implies $g(x) <0$, 
we have $g( \tr P) < 0$.  Thus $ \l_{\rm min}(P) \geq 0 > g( \tr P)$
which implies that $P\in \Int M^g$, since $g$ is continuous. \qed

\Cor{\LL.3}  {\sl
The dual subequation $\wt{M^g}$ is borderline.
}
\pf
The first part of (\LL.2) implies that $0\in \partial M^g = - \partial \wt { M^g}$.
Combined  with the second  part of (\LL.2), this is condition (1)$'$ in Lemma \CC.2
for the subequation $F=\wt{M^g}$, which proves that $\wt{M^g}$ is borderline.\qed

\Prop{\LL.4}  {\sl
The subequation $M^g$ is additive, i.e., $M^g+M^g \ss M^g$, if and only if
$g$ is subadditive, i.e., $g(x+y) \leq g(x) + g(y)$.
}
\pf
Use (\LL.3) and $\tr(A+B) = \tr A + \tr B$.\qed
\medskip

If $g(x) \equiv - \d x \ (\d>0)$, then $M^g  \equiv  \cp(\d)$ is the convex cone
subequation discussed in Appendix B.  However, there are plenty of other subadditive 
decreasing functions $g$.

Suppose $g$ is concave on $[0,a]$ with $g(0)=0$.  Then, as noted in the introduction to
[\BRUCK], the extension of $g(x)$ from $[0,a]$ to $[0,\infty)$ defined by 
$$
g(x) \ \equiv\ j g(a) + g(x-ja), \qquad ja\leq x\leq (j+1)a, \qquad j=1,2,...
\eqno{(\LL.4)}
$$
is subadditive on $[0,\infty)$ and has the property that $g\geq h$ for
any other subadditive function $h$ on $[0,\infty)$ which agrees with $g$
on $[0,a]$.  The elementary proof is omitted. Summarizing, we have the following.

\Theorem{\LL.5} {\sl
Suppose that $g:[0,\infty) \to \bbr$ is the extension of a decreasing concave
function on $[0,a]$ defined by (\LL.4) with $g(0)=0$.  Then $M^g$
is a monotonicity subequation (orthogonally invariant), and its dual 
$\wt{M^g}$ is borderline.
}
 
\Lemma {\LL.6}  {\sl
The dual subequation $\wt{M^g}$ is defined by
$$
\wt{M^g}\ \ :\quad \tr A \geq 0 \ \ {\rm or\ \ }  \l_{\rm max}(A) \ \geq \ - g\left ( -  \tr A\right) \quad( {\rm with}\ \ \tr A \leq0)
$$
}
\pf Note that
$$
\eqalign
{
A\in \wt{M^g}  \quad \iff \quad -A \notin \Int M^g 
  & \quad \iff \quad  \l_{\rm min}(-A) \ \leq \  g\left ( -  \tr A\right) 
 \quad {\rm or}\quad \tr(-A)\leq 0
 \cr
 & \quad \iff \quad \l_{\rm max}(A) \ \geq \ - g\left ( -  \tr A\right) 
  \quad {\rm or}\quad \tr(A)\geq 0
}
$$
since $ \l_{\rm max}(A) = -  \l_{\rm min}(-A)$.  \qed

\Prop {\LL.7} {\sl
The characteristic function $f$ for the dual subequation  $\wt{M^g}$ on $\rn$  is 
$
f(\l) =  g^{-1} (-\l) + (n-1)\l
$
for $\l\geq0$.
}
\pf
The increasing radial profile of $\wt{M^g}$ is by definition
$$
\L 
 \equiv\ \{(\l, \mu) : \l P_{e^\perp} + \mu P_e \in \wt{M^g} \ \ {\rm and}\ \ \l\geq0\}.
$$
Note that $\tr A = (n-1)\l +\mu$
if $A \equiv \l P_{e^\perp} + \mu P_e$.
If $\l \geq0$ and $A \in  \wt{M^g}$ with $\tr A \leq0$, then
$$
\l \ \equiv \ \l_{\rm max} \ \geq\ 0, \quad \mu\ \leq\ 0, \qquad 
{\rm and\  hence}\qquad   \l \geq -g\left(- (n-1)\l - \mu\right)
$$
Set $x\equiv - (n-1)\l - \mu\geq0$
and 
$y \equiv -\l \leq0$.
Then $y\leq g(x)$ is equivalent to $x\leq g^{-1}(y)$ since $g$ is decreasing and $g(0)=0$.
Thus $ -(n-1)\l -\mu \leq g^{-1}(-\l)$, or 
$\mu +  g^{-1}(-\l) + (n-1) \l \geq0$.
Since $f$ is defined by $\mu+f(\l)\geq0$ for such pairs $(\l,\mu)$,this completes the proof.\qed

\Ex{\LL.8. (An Explicit Example where (SC) holds but the Subequation is not Contained in a Uniformly Elliptic Subequation)}  Define $g:[0,a]\to [-b,0]$ via its inverse by
$$
g^{-1}(-\l) \ \equiv\ \l(\a-2\log \l) \qquad 0\leq\l\leq a.
\eqno{(\LL.5)}
$$
Here $\a$ is a constant chosen first, and then $a$ is chosen small enough so that 
$h(\l) \equiv g^{-1}(-\l)$ is strictly increasing on $[0,a]$, and finally we set $-b = g(a)$.
Note that $h'(\l) = \a -2-2\log \l$.
Also, $h''(\l) = -{2\over \l} <0$. Therefore $g$ is concave and strictly decreasing
on $[0,a]$ with $g(0)=0$. Applying Theorem \LL.4 we see that 
$$
M^g \ \ {\rm is \ a\  monotonicity\  subequation\  whose\  dual\  } \wt{M^g}\ {\rm is\ borderline.}
\eqno{(\LL.6)}
$$
By Proposition \LL.7
$$
{\rm The \  dual\  } \wt{M^g}\ {\rm has\  characteristic\ function\ } f(\l) = \l  (\a + n-1 -2\log \l) \ {\rm on}\ [0,a].
\eqno{(\LL.7)}
$$

Recall the subequation $H(\a')$ discussed in Example \JJ.11.  If we take $y=\l$, $k=2$, and the $\a$ there 
to be the $\a' \equiv \a+n-1$ for (\LL.7), then the characteristic function $f(\l)$ for $\wt{M^g}$
is the same as the characteristic function for $H(\a')$, for $\l$ small.  Since $\int_{0^+} {1\over f} =\infty$,
as shown there, this proves
\Prop{\LL.9} 
$$
{\sl The \  (SMP) \ holds\ for\  this\ dual\ subequation\ } \wt{M^g},\ {\sl and \ so\   the \  (SC) \ holds\ for\  } M^g.
%\eqno{(\LL.9)}
$$

It is easy to see that $\wt{M^g}$ is not contained in a uniformly elliptic subequation
since $f(\l)/\l = \a +n-1 -2\log \l \to \infty$ as $\l\to\infty$.

Finally we remark that, as in Example \JJ.11, if $\b\equiv e^{{1\over 2}(1+\a)}$, then
the Hopf function
$$
\psi(|x|) \ \equiv\    e^{-\b R^2/2} -  e^{-\b |x|^2/2} \ \ {\rm is}\ \wt{M^g}\,{\rm harmonic \ for\ } |x|\ {\rm small}.
\eqno{(\LL.8)}
$$
 (This function $\psi(|x|)$ is also a harmonic for the subequations
$F_f^{\rm min/2} \ss \wt{M^g} \ss F_f^{\rm min/max}$ described in Theorem \JJ.5.)

 \vfill\eject
%\vskip.3in

%%%%%%%%%%%%%%%%%%%%%%%%%%%%%%%%%%%%%%%%%% 
%%%%%%%%%%%%%%%%%%%%%%%%%%%%%%%%%%%%%%%%%% 
%%%%%%%%%%%%%%%%%%%%%%%%%%%%%%%%%%%%%%%%%% 
%%%%%%%%%%%%%%%%%%%%%%%%%%%%%%%%%%%%%%%%%% 
%%%%%%%%%%%%%%%%%%%%%%%%%%%%%%%%%%%%%%%%%% 

\noindent{\headfont \MM.\  Another Application  --  Product Subequations}.
\medskip
In this section we apply our main result to study the (SMP) for product subequations.
Let $F\ss\Symn$ and $G\ss \Sym(\bbr^m)$ be invariant pure second-order subequations,
%on $\rn$ and $\bbf^m$ respectively. 
and  consider the product subequation  $H \equiv $  ``$F\times G$'' $\ss \Sym(\bbr^{n+m})$ 
defined, for coordinates $(x,y) \in \rn\times \bbr^m$, by requiring that 
$D^2_x u\in F$ and $D^2_y u\in G$. In other words, $u$ is separately $F$-subharmonic
in $x$ and $G$-subharmonic in $y$.

One easily checks that if either $F$ or $G$ is stable, then $H$ is stable. On the other hand,
 if either $F$ or $G$ has a counterexample to the (SMP), so does $H$ (take the same counterexample
 considered as a function of all the variables).
For the remaining case we have the following. The proofs are omitted.

 \Theorem{\MM.1}  {\sl
 Let $F$ and $G$ be invariant borderline subequations for which the (SMP) holds.
Let $f$ and $g$ denote their   respective characteristic functions.
 Suppose one of these, say $g$, satisfies
 $$
 {g(y)-g(x) \over  y-x} \ \ \ {\rm is\ bounded\  for }\ \ 0<x<y\ \ {\rm small}.
 \eqno{(\MM.1)}
 $$
Then the  (SMP) holds  for the product  subequation $H=F\times G$.}
\medskip
\noindent
{\bf Outline of Proof.}
The first step is the following.
 
 \Prop{\MM.2}  {\sl Let $\overline h$ and $\underline h$
be the upper and lower characteristic functions of $H$.
 Then}
 $$
 \overline h (\l) \ =\ f(\l)+g(\l) + \l
 \and
 \underline h(\l)  \ =\ \min\{f(\l),g(\l)\}
 $$
 \noindent
The next step is that
 $$
 \eqalign
 {
 \int_{0^+} {1\over f} \ =\ \infty
 \quad{\rm and}\quad 
  \int_{0^+} {1\over g}  \ =\ \infty 
   \quad{\rm and}\quad (\MM.1)
 \qquad&\Rightarrow\qquad  \int_{0^+} {1\over f+g} \ =\ \infty   \cr
 {\rm and\ therefore} 
  \qquad&\Rightarrow\qquad  \int_{0^+} {1\over f+g+\l} \ =\ \infty  
  }
 $$
 Theorem \MM.1 now follows from Theorem A$'$. \qed
 
 \Remark{\MM.3}  There  exist functions $f, g>0$ with 
 $\int_{0^+} {1\over f} = \int_{0^+} {1\over g} = \infty$ but $\int_{0^+} {1\over f+g} < \infty$.

 \vfill\eject
%\vskip.3in

%%%%%%%%%%%%%%%%%%%%%%%%%%%%%%%%%%%%%%%%%% 
%%%%%%%%%%%%%%%%%%%%%%%%%%%%%%%%%%%%%%%%%% 
%%%%%%%%%%%%%%%%%%%%%%%%%%%%%%%%%%%%%%%%%% 
%%%%%%%%%%%%%%%%%%%%%%%%%%%%%%%%%%%%%%%%%% 
%%%%%%%%%%%%%%%%%%%%%%%%%%%%%%%%%%%%%%%%%% 

\centerline{\headfont Appendix A.  Radial Subharmonics}.
 \smallskip

Since our characterization of radial subharmonics is useful for many purposes,
it is separated out in this appendix.  Recall  the characteristic lower function $\uf$ associated 
with a subequation $F$ and  the radial subequation $R_{\uf}$ defined by
$$
\psi'' + \uf \left({\psi'  \over t}\right)\ \geq\ 0\qquad{\rm on} \ 0<t<\infty.
$$
In the following we drop the bar, letting $f$ denote $\uf$.

\Theorem{A.1. (Radial Subharmonics)} {\sl
The function $u(x) \equiv \psi(|x|)$ is $F$-subharmonic on an  annular region in $\rn$
if and only if 
 $\psi(t)$ is $R_f$-subharmonic on the corresponding sub-interval of $(0,\infty)$.
 }
\medskip
\noindent
{\bf Proof.  ($\Rightarrow$):}  Suppose $u(x) \equiv \psi(|x|)$ is $F$-subharmonic.
If $\vf(t)$ is a test function for $\psi(t)$ at $t_0$, then $\vf(|x|)$ is a test function
for  $\psi(|x|)$ at any point on the $t_0$-sphere in $\rn$.
Therefore $D^2_{x_0}\vf \in F$.  Applying the formula  (Lemma \EE.1)  for $D^2_{x_0}\vf$ in terms of 
$\vf'(t_0)$ and $\vf''(t_0)$, the equivalence (\EE.1),
 and the definition of $(R_F)_{t_0}$,  we have $J^2_{t_0} \vf(t) \in R_F$.
This proves that $\psi(t)$ is $R_F$-subharmonic. 

\noindent
{\bf  ($\Leftarrow$):}   Suppose that $\psi(t)$ is $R_F$-subharmonic.  We must show 
that $u(x) \equiv \psi(|x|)$ is $F$-subharmonic. That is, given a test function $\vf(x)$ for
$u(x)$ at a point $x_0$, we must show that $D^2_{x_0}\vf \in F$.

Suppose that there exists a {\sl smooth} function $\overline\psi(t)$, defined near $t_0=|x_0|$,
such that $\overline\vf(x) \equiv \overline\psi(|x|)$ satisfies
$$
u(x)\ \leq\ \overline\vf(x) \ \leq \vf(x)
\eqno{(A.1)}
$$
near $x_0$.  Then $\overline\psi(t)$ is a test function for $\psi(t)$ at $t_0$.
Hence, the 2-jet of $\overline\psi$ at $t_0$ belongs to $R_F$. 
By Lemma \EE.1 and the discussion above, this implies that $D^2_{x_0}\overline\vf \in F$.
The inequality $\overline\vf(x) \leq\vf(x)$ (with equality at $x_0$) implies
that $D^2_{x_0}\vf = D^2_{x_0}\overline\vf +P$ for some $P\geq0$,
which proves that $D^2_{x_0}\vf\in F$ as desired.

To complete this argument by finding $\overline\psi(t)$ there is some flexibility given by Lemma 2.4 in [\DDR]
so that not all test functions $\vf(x)$ need  be considered.  
First we may choose new   coordinates
$z=(t,y)$ near $x_0$ so that $t\equiv |x| $. (Thus $t=$   constant defines the sphere of radius $t$ near $x_0$.)
Furthermore, we may assume that $\vf(z)$ is a polynomial of degree $\leq 2$ in 
$z=(t,y)$ and that it is a {\sl strict} local test function, i.e., $u(z)  < \vf(z)$ for $z\neq z_0$.
Now Lemma A.2 below  ensures the existence of $\overline\vf(x) = \overline\psi(|x|)$
satisfying (A.1). \qed

\medskip

Let $z=(t,y)$ denote standard coordinates on $\rn=\bbr^k\times \bbr^\ell$.
Fix a point $z_0=(t_0,y_0)$ and let $u(t)$ be an upper semi-continuous function
(of $t$ alone) and $\vf(z)$ a $C^2$-function, both defined in a neighborhood of $z_0$.

\Lemma {A.2} {\sl
Suppose $u(t) < \vf(z)$ for $z\neq z_0$ with equality at $z_0$. 
If $\vf(z)$ is a polynomial of degree $\leq2$, then there exists 
a polynomial $\overline\vf(t)$ of degree $\leq2$ with
$$
u(t)  \ \leq\ \overline\vf(t) \ \leq\ \vf(z) 
\qquad{\rm near}\ \ z_0.
\eqno{(A.2)}
$$
}
\pf
We may assume $z_0=0$ and $u(0)=\vf(0)=0$. Then
$$
\vf(z) \ =\  \bra pt + \bra qy + \bra {At} t + 2 \bra {Bt} y + \bra {Cy}y.
$$
We assume $u(t)  < \vf(t,y)$ for $|t| \leq \e$ and $|y| \leq \d$ with $(t,y)\neq (0,0)$.

Setting $t=0$,  we have $0=u(0) < \bra qy +\bra {Cy}y$
for $y\neq 0$ sufficiently small.  Therefore, $q=0$ and $C>0$ (positive definite). 
Now define
$$
\overline\vf(t)\ \equiv\ \bra pt + \bra {(A-B^tC^{-1} B)t}{t}.
\eqno{(A.3)}
$$
The inequalities in (A.2) follow from the fact that for $t$ sufficiently small,
$$
\overline\vf(t)\ = \ \inf_{|y|\leq \d} \vf(z)\ =\ \bra  pt + \bra {At}t + \inf_{|y|\leq \d} \{ 2\bra {Bt} y+\bra{Cy}y\}.
\eqno{(A.4)}
$$
To prove (A.4) fix $t$ and consider the function $2\bra {Bt} y+\bra{Cy}y$.
Since $C>0$, it has a unique minimum point at the critical point $y= - C^{-1}Bt$. 
The minimum value is $-\bra{B^tC^{-1}Bt} t$. If  $t$ is sufficiently small, the critical point
$y$ satisfies $|y|<\d$, which proves (A.4).\qed

 %\vfill\eject
\vskip.3in

%%%%%%%%%%%%%%%%%%%%%%%%%%%%%%%%%%%%%%%%%% 
%%%%%%%%%%%%%%%%%%%%%%%%%%%%%%%%%%%%%%%%%% 
%%%%%%%%%%%%%%%%%%%%%%%%%%%%%%%%%%%%%%%%%% 
%%%%%%%%%%%%%%%%%%%%%%%%%%%%%%%%%%%%%%%%%% 
%%%%%%%%%%%%%%%%%%%%%%%%%%%%%%%%%%%%%%%%%% 

\centerline{\headfont Appendix B.  Uniform Ellipticity}.
 \medskip

This is a geometric  discussion of uniform ellipticity.
A family of convex cone subequations $\{M_\d\}$ is said to be a 
{\sl fundamental neighborhood system} for $\cp$ if given any conical
neighborhood $G$ of $\cp$ (this means that $\cp-\{0\} \ss \Int G$ and $G$ is a cone),
there exists $\d$ with $M_\d\ss G$.
 Given such a family $\{M_\d\}$, 
a subequation $F$ is {\bf uniformly elliptic} if one of the  $M_\d$ is a monotonicity 
subequation for $F$.  That is,
$$
F+M_\d\ \ss\ F  \qquad{\rm for\ some\ \ } \d.
\eqno{(B.1)}
$$
This definition is easily seen to be  independent  of the choice of the neighborhood system
$\{M_\d\}$  for $\cp$.
(The monotonicity condition (B.1) can always be rephrased classically,
in terms of the operator defining $M_\d$, as two
inequalities -- see, for example, (4.5.1)$'$ in [\RS]).

The standard choice made in the literature consists of the {\bf Pucci cones}
\smallskip
\centerline
{
$
\cp_{\l,\L} \ \equiv\ \{A : \l\tr\, A^+ + \L\tr\, A^- \geq0\}
$
}\smallskip\noindent
with $0<\l<\L$, where $A=A^++A^-$ is the decomposition of $A$ into positive and negative parts.
Another good choice is the {\bf  $\d$-uniformly elliptic regularization} $\cp(\d)$ of $\cp$
$$
\cp(\d) \ \equiv\ \{A : A+\d  ( \tr\, A)I\geq0\} \qquad(\d>0).
$$
Both $\cp_{\l,\L}$ and $\cp(\d)$ are convex cone subequations as required.
See Section 4.5 of [\SURVEY] for more details regarding $\cp_{\l,\L}$ and $\cp(\d)$
(The Riesz characteristics are computed in Example 6.2.5.)

Since there is a largest monotonicity subequation $M_F$ for $F$, uniform ellipticity can be
defined equivalently as 
$$
M_F \ \ {\rm contains\ a\ convex\ conical \  neighborhood \ of \ } \cp,\ \ {\rm or\ as}
\eqno{(B.1)'}
$$
$$
\wt M_F \ \ss\ \wt{\cp(\d)}\ \ \ {\rm for \ some\ \ }\d>0.
\eqno{(B.1)''}
$$
 
\Lemma {B.1} {\sl
If $F$ is normal, i.e., $M_F$ is a cone (and hence a convex cone), then\medskip

\centerline{
$F$ is uniformly elliptic $\iff$ $M_F$ is a conical neighborhood of $\cp \iff$
}
\medskip
\centerline{  $P_e\in\Int M_F$
for all $e\neq0$ $\iff -P_e\notin \wt M_F$ for all $e\neq0 \iff \wt M_F$ is borderline.}
}
\medskip
The next remark is to be applied to $F=M_G$ where $G$ is normal.

\Remark{B.2. (Cone subequations and the Riesz Characteristic)} For simplicity suppose that $f=\uf=\of$ is the characteristic function for a cone 
subequation $F$.  Then $f(t\l) = tf(\l)$ for $t>0$, and hence the characteristic  function reduces
to two numerical invariants
$$
\a \ \equiv\ f(1)  \and    \a^*\ \equiv\ -f(-1), \qquad 0\leq \a, \a^* \leq \infty
\eqno{(B.2)}
$$
where we have
$$
f(\l) \ =\ \a \l\quad {\rm for}\ \l>0 \and  f(\l) \ =\   \a^* \l \quad {\rm for}\ \l<0.
\eqno{(B.3)}
$$
The radial profile $\L$  is defined by
$$
\mu + \a \l \ \geq\ 0\ \ \ {\rm if}\ \l\geq0 \and \mu + \a^*\l \ \geq\ 0 \ \ {\rm if}\ \l\leq0.
\eqno{(B.4)}
$$
Note that  $\a=\infty \iff P_{e^\perp} -\mu P_e \in F$ for all $\mu \iff -P_e \in F \iff F$ is
not borderline.  That is,
$$
F\ \ {\rm satisfies\ the \ (SMP)} \qquad\iff\qquad \a \equiv \a_F\ <\ \infty.
\eqno{(B.5)}
$$
The invariant $p_F \equiv \a_F+1$ is called the {\bf Riesz characteristic of $F$}
because of its connection with Riesz kernels.  See [\RS], [\ASPECTS] for applications, 
examples and a fuller discussion, where it is proved, in particular, that 
$\a \a^*\geq1$.

% \vfill\eject
\vskip.3in

%%%%%%%%%%%%%%%%%%%%%%%%%%%%%%%%%%%%%%%%%% 
%%%%%%%%%%%%%%%%%%%%%%%%%%%%%%%%%%%%%%%%%% 
%%%%%%%%%%%%%%%%%%%%%%%%%%%%%%%%%%%%%%%%%% 
%%%%%%%%%%%%%%%%%%%%%%%%%%%%%%%%%%%%%%%%%% 
%%%%%%%%%%%%%%%%%%%%%%%%%%%%%%%%%%%%%%%%%% 

\centerline{\bf References}

\vskip .2in

\noindent
\item{[\AS]}  S. N.   Armstrong, B.  Sirakov  and  C. K. Smart,  {\sl Fundamental
solutions of homogeneous fully nonlinear elliptic equations}, Comm.
Pure. Appl. Math.  {\bf 64} (2011), 737-777.

\smallskip

 \item{[\BD]}  M. Bardi and F. Da Lio,  {\sl
On the strong maximum principle for fully nonlinear degenerate elliptic equations}, 
Arch. Math. {\bf 73} (1999), 276-285.

 \smallskip

 \item{[\BaB]}  G. Barles and J. Busca,  {\sl
Existence and comparison results for fully nonlinear degenerate elliptic equations
without zeroth-order term},   Comm. in Partial Diff. Eqs.
 {\bf 26} (2001), 2323-2337.

 \smallskip

 \item{[\BRUCK]}  A. Bruckner,  {\sl
Minimal superadditive extensions of superadditive functions}, 
Pacific J. Math. {\bf 10} (1960), 1155-1162.

 \smallskip
 
\item {[\CaC]}  L.  Caffarelli and X. Cabr\'e,  {\sl  Fully nonlinear elliptic equations},   
Colloquium Publications {\bf 43}, American Math. Soc., Providence, 1995.

 \smallskip
 
\item {[\CLN]}  L.  Caffarelli,  Y.Y.  Li and L. Nirenberg,  {\sl  Some remarks on singular solutions
of nonlinear elliptic equations, III: viscosity solutions,
including parabolic operators},  Comm. Pure Appl. Math. {\bf 66} (2013), 109Ð143. 
 ArXiv:1101.2833.

\smallskip

\noindent
\item{[\C]}   M. G. Crandall,  {\sl  Viscosity solutions: a primer},  
pp. 1-43 in ``Viscosity Solutions and Applications''  Ed.'s Dolcetta and Lions, 
SLNM {\bf 1660}, Springer Press, New York, 1997.

 \smallskip

\noindent
\item{[\CIL]}   M. G. Crandall, H. Ishii and P. L. Lions {\sl
User's guide to viscosity solutions of second order partial differential equations},  
Bull. Amer. Math. Soc. (N.S.) {\bf 27} (1992), 1-67.

 \smallskip

\noindent
\item{[\EB]}   A. Eberhard, {\sl
 Prox-regularity and subjects},  in A. Rubinov (Ed.), ``Optimization and Related Topics'',
 Applied Optimization Volumes, Kluwer Academic Publishers, Dordrecht, 2001,  pp. 237-313.

 \smallskip

 \item{[\F]} H. Federer, Geometric Measure Theory, Springer-Verlag, New York, 1969.
 
 \smallskip

%   \noindent
%\item{[\K]}    N. V. Krylov,    {\sl  On the general notion of fully nonlinear second-order elliptic equations},    %Trans. Amer. Math. Soc. (3)
% {\bf  347}  (1979), 30-34.
%\smallskip

\item {[\DDD]}   F. R. Harvey and H. B. Lawson, Jr., {\sl  Dirichlet duality and the non-linear Dirichlet problem},    Comm. on Pure and Applied Math. {\bf 62} (2009), 396-443. 

ArXiv:math.0710.3991.

\smallskip

\item {[\DDR]}  \ \----------, {\sl Dirichlet Duality and the Nonlinear Dirichlet Problem on Riemannian Manifolds},  J. Diff. Geom. {\bf 88} (2011), 395-482.   ArXiv:0907.1981.

\smallskip

%\item  {[\HYP]} \ \----------, {\sl  Hyperbolic polynomials and the Dirichlet problem},   ArXiv:0912.5220.

%\smallskip

%\item {[\HLGAR]}  \ \----------, {\sl  G\aa rding's theory of hyperbolic polynomials},
%   {Communications in Pure and Applied Mathematics} (to appear).  

% \smallskip

 \item{[\SURVEY]}  \ \----------, 
 {\sl  Existence, uniqueness and removable singularities
for nonlinear partial differential equations in geometry},\ 
(with  R. Harvey),  fpp. 102-156 in ``Surveys in Geometry 2013'', vol. 18,  
H.-D. Cao and S.-T. Yau eds., International Press, Somerville, MA, 2013.
ArXiv:1303.1117.

\smallskip

\item {[\RS]}  \ \----------,  {\sl  Removable singularities for nonlinear subequations},  \  Indiana Univ. Math. J., {\bf 63}, No. 5 (2014), 1525-1552.
 ArXiv:1303.0437.

\smallskip

\item  {[\ASPECTS]} \ \----------,      {\sl Tangents to subsolutions -- existence and uniqueness, Parts I and II}, 

ArXiv:1408.5797 and ArXiv:1408.5851.

\smallskip

\item  {[\AETHM]} \ \----------,  {\sl  The AE Theorem and addition theorems for quasi-convex functions},  

ArXiv:1309.1770.

 \smallskip

   \noindent
\item{[\HIR]}    J.-B. Hiriart-Urruty and  Ph. Plazanet,    {\sl  Moreau's Theorem revisited},  
 Analyse Nonlin\'eaire (Perpignan, 1987), Ann. Inst. H. Poincar\'e {\bf 6} (1989), 325-338.

\smallskip

   \noindent
\item{[\Kawohl]}    B. Kawohl and N. Kutev,    {\sl  Strong maximum principle for semicontinuous 
viscosity solutions of nonlinear partial differential equations},    Arch. Math.
 {\bf  70}  (1998), 470-478.

\smallskip

%   \noindent
%\item{[\Kawohll]}  \ \----------,   {\sl  Comparison principle for viscosity solutions of
%fully nonlinear, degenerate elliptic equations},    Comm. in Partial Diff. Eqns.
 %{\bf  32}  (2007), 1209-1224.

%\smallskip

   \noindent
\item{[\Krylov]}    N. V. Krylov,    {\sl  On the general notion of fully nonlinear second-order elliptic equations},    Trans. Amer. Math. Soc. (3)
 {\bf  347}  (1979), 30-34.

\smallskip

   \noindent
\item{[\LE]}    Y.  Luo  and A. Eberhard,    {\sl Comparison principle for viscosity solutions of elliptic equations
via fuzzy sum rule},    J. Math. Anal. Appl. 
 {\bf  307}  (2005), 736-752.

\smallskip

\end